\documentclass[review,3p,times,numbers,sort&compress]{elsarticle}

\usepackage{amssymb}
\usepackage{amsmath}
\usepackage{amsfonts}
\usepackage{bm}
\usepackage{amsthm}
\newtheorem{theorem}{Theorem}
\newtheorem{proposition}{Proposition} 

\usepackage[ruled,vlined,linesnumbered]{algorithm2e}
\usepackage{multirow}
\SetKwInput{KwIn}{Require}
\SetKwInput{KwOut}{Ensure}

\usepackage{caption}
\usepackage{subcaption}
\usepackage{xcolor}    
\usepackage{hyperref}  
\usepackage{fontawesome5}
\usepackage{enumitem}
\usepackage{float}
\usepackage{placeins}
\usepackage{titlesec}

\titleformat{\paragraph}[runin]
{\itshape}
{}
{0pt}
{}

\begin{document}

\begin{frontmatter}

\title{GeoFunFlow-3D: A Physics-Guided Generative Flow Matching Framework for High-Fidelity 3D Aerodynamic Inference over Complex Geometries}
\author[a]{Ruiling Jiang} 
\author[a]{Yong Zhang\corref{cor1}} 
\author[a]{Houbiao Li\corref{cor2}}

\cortext[cor1]{First corresponding author. \newline E-mail address: mathzy@163.com}  
\cortext[cor2]{Second corresponding author. \newline E-mail address: lihoubiao0189@163.com}

\address[a]{School of Mathematical Sciences, University of Electronic Science and Technology of China, Chengdu, 611731, China}

\begin{abstract}
Deep generative models and neural operators have demonstrated significant potential for 3D aerodynamic inference. However, they often face inherent challenges in maintaining physical consistency and preserving high-frequency features, primarily due to spectral bias and gradient conflicts within the governing equations. To address these issues, we propose GeoFunFlow-3D, a physics-guided generative flow matching framework. Temporally, we utilize optimal transport theory to build the generation path, ensuring stable training dynamics. Spectrally, we introduce a high-order discrete engine without automatic differentiation (No-AD) to reduce gradient stiffness. Spatially, a topology-aware super-resolution module (SATO) is employed to rigorously enforce physical laws in localized regions such as shock waves. We evaluated our framework on complex industrial datasets. On the BlendedNet dataset, the model successfully avoids mode collapse even under sparse data conditions. For the NASA Rotor37 test, it accurately captures 3D detached shock structures. Compared to conventional operators, GeoFunFlow-3D significantly improves accuracy, reducing the pressure field error (RRMSE) to 0.0215 while maintaining competitive inference efficiency. Ultimately, this work provides a reliable, geometry-driven approach for generating high-dimensional fluid fields. 
\end{abstract}


\end{frontmatter}

\section{Introduction}

\subsection{Problem Statement}
Extending physics-informed neural networks (PINNs) and neural operators to 3D aerodynamic forward generation tasks remains a challenging task. Existing models still encounter specific bottlenecks in terms of computational rigor and physical consistency:

\begin{enumerate}
    \item \textbf{Memory Complexity and Spectral Bias:} While classical continuous models have achieved foundational success, their reliance on automatic differentiation (AD) for high-order derivatives often results in substantial memory footprint overhead in 3D high-resolution scenarios. Mathematically, continuous AD mechanisms induce spectral distortion based on the Neural Tangent Kernel (NTK) theory~\cite{wang2022when}. This can contribute to Spectral Bias for high-frequency physical features inside boundary layers.
    \item \textbf{Boundary Constraint Challenges and Gradient Conflict:} Global partial differential equation (PDE) residuals often conflict with local boundary conditions. Such inconsistencies are especially pronounced in regions with complex geometries, such as sharp edges or compressor interiors, and may lead to non-physical artifacts, including cross-boundary mass leakage. Thus, the thermodynamic conservation laws of transonic internal flow fields remain challenging to model rigorously.
    \item \textbf{Dilemma of Zero-Observation Forward Generation:} Flow matching models rely heavily on large amounts of high-fidelity data to build probability measures~\cite{wang2025geofunflow}. Physical observations are lacking in early aerospace design. Thus, PDEs must act as physical penalties to replace missing data for pure geometry-prior generation. Unlike traditional unsupervised learning or purely physics-driven solvers, zero-observation forward generation establishes a probability distribution based strictly on geometric priors via generative flow matching, without relying on any downstream task annotations. However, without proper scheduling, rigid nonlinear physical constraints can disrupt the originally smooth data probability manifold. This causes strong gradient conflicts and local minima during optimization. Consequently, such ill-conditioning of the PDE system and numerical stiffness can inadvertently introduce non-physical artifacts.
\end{enumerate}

\subsection{Primary Contributions}
To resolve these limitations, we propose GeoFunFlow-3D, a physics-guided 3D fluid generation foundation model. Our previous work, GeoFunFlow~\cite{wang2025geofunflow}, established a geometric flow matching framework for inverse problems, which inherently relies on sparse sensor observations to reconstruct physical fields. In contrast, this study targets the more challenging zero-observation forward generation for 3D aerodynamics, where the inference phase relies exclusively on geometric priors without any target physical observations. We directly address the severe dimensionality curse and rigid physical boundary constraints. Our core contribution lies in the construction of a unified spectral-spatial-temporal convergence framework. This mechanism achieves a theoretical leap from observation-dependent posterior sampling to pure geometry-prior forward generation:

\begin{enumerate}
    \item \textbf{Spectral Stability: Linear stability guarantee based on the No-AD engine.} We introduce a high-order discrete difference engine without automatic differentiation (No-AD) and mathematically prove the boundedness of its NTK eigenvalues (see Theorem 1 and its proof in Section 4.3). This design effectively mitigates Hessian rigidity and alleviates high-frequency spectral bias.
    \item \textbf{Spatial Consistency: Theoretical alignment of SATO and neural MsFEM.} We propose a topology-aware super-resolution module (SATO). Its residual compensation mechanism aligns with the over-sampling theory of the multiscale finite element method (MsFEM)~\cite{hou1997multiscale}. Supported by a phase-field hard mask, it ensures rigorous spatial local consistency and strictly enforces the non-penetration condition within complex boundary layers.
    \item \textbf{Global Convergence: Path smoothing based on Variational Homotopy Scheduling and exponential relaxation.} Following the numerical continuation method~\cite{allgower2012numerical}, flow matching is built on an approximate $2$-Wasserstein ($W_2$) geodesic~\cite{benamou2000computational}. We derive a $C^1$ continuous Variational Homotopy Scheduling protocol and an exponential physical relaxation mechanism. This ensures the energy optimality and global nonlinear Well-posedness of the probability distribution during optimization.
\end{enumerate}

\section{Related Work}

GeoFunFlow-3D builds upon neural operators, physics-informed learning, Lagrangian boundary awareness, and generative flow matching. This section reviews these technologies and analyzes how the architecture is rebuilt for 3D aerodynamic inference.

\subsection{Neural Operators and Geometric Feature Enhancement}
Neural operators construct mappings between infinite-dimensional spaces. Li et al.~\cite{li2021fourier} proposed the Fourier Neural Operator (FNO) to decouple mappings from resolutions. To handle unstructured point clouds, Li et al.~\cite{li2023geometry} proposed the Geometry-Informed Neural Operator (GINO). Furthermore, the Feature-Enhanced Neural Network (FENN) proposed by Song et al. ~\cite{song2025fenn} demonstrated that, rather than having the network learn geometry from scratch, explicitly encoding geometric metrics such as the symbolic distance function (SDF) from a point to a surface—as prior knowledge significantly reduces the difficulty of optimizing physical equations. The GP-DE-PINN framework proposed by Wang et al.~\cite{wang2026gppinn}  further validated the effectiveness of decoupling geometric encoding from the evolution of physical fields.

\textbf{Improvements in this study:} GeoFunFlow-3D combines GINO and FNO (Sections 3.1 and 3.4). Feature injection is expanded into a 9D multimodal geometric embedding (including normal vectors, curvatures, and SDF). This assigns complete Riemannian manifold tangent bundle properties to point clouds. We construct a joint mechanism of front-end homeomorphic dimension reduction and back-end continuous decoding. This improves topological generalization for complex flow fields.

\subsection{Physics-Informed Deep Learning and Complex Manifold Mapping}
Traditional convolutional neural networks (CNNs) can introduce non-conformal errors on irregular 3D topologies. Gao et al.~\cite{gao2021phygeonet} proposed PhyGeoNet using isogeometric analysis (IGA). By solving elliptic partial differential equations, a coordinate mapping is constructed from the physical domain to a regular reference domain, thereby enabling the use of efficient difference convolution as an alternative to the extremely memory-intensive automatic differentiation (AD). De Ryck et al.~\cite{ryck2022wpinns} demonstrated in their work on weak-form PINNs (wPINNs)  that, when dealing with sharp discontinuities such as shock waves, the mathematical significance of the discrete weak-form is more important than the derivatives of the pointwise strong-form. While PhyGeoNet elegantly mitigates AD overhead, adapting its explicit coordinate transformations to 3D topologies with extreme curvature changes (such as wing-body fuselages) may encounter challenges related to local grid distortions and Jacobian matrix singularities.

\textbf{Improvements in this study:} We derive a No-AD discrete differential engine in the latent space (Section 3.3) and pull back gradients using the inverse Jacobian matrix. Instead of relying on explicit global mesh transformations, this operator acts inherently as a band-limited filter (the theoretical proof regarding how it suppresses the $\mathcal{O}(\omega^{2k})$ NTK eigenvalue explosion and eliminates extreme Hessian Rigidity is detailed in Theorem 1, Section 4.3).

\subsection{Lagrangian Topology Awareness and Absolute Boundary Hard-Constraining}
Classical PINNs use density-based soft penalties. This can occasionally lead to blurred boundaries or inadequate enforcement of Dirichlet conditions. Zhou et al.~\cite{zhou2026ltpinn} proposed LT-PINN to integrate PINNs with constructive solid geometry (CSG). LT-PINN uses an SDF mask to strictly confine the PDE solution to the fluid domain and precisely extracts the unit outward normal vector via an analytical function, thereby fundamentally eliminating the “directional ambiguity” of physical flux constraints. Although LT-PINN has overcome the theoretical barrier between pseudo-boundaries and analytical boundaries, its parameterization process, which relies on the combination of multiple elements, can become computationally complex to scale when applied to industrial three-dimensional aircraft models consisting of millions of patches.

\textbf{Improvements in this study:} We upgrade the topology-aware core (Section 3.4). By introducing a phase-field hard mask and tangent bundle projection, we implement Euler wall tangency constraints. This mathematically avoids the cross-boundary leakage shown in the Ablation Study. Additionally, the SATO module is proposed as a neural generalization of MsFEM. It dynamically generates microscopic basis functions from high-frequency features, eliminating resonance errors between grid scales.

\subsection{Generative Flow Matching and Dual Optimization Dynamics}
Deep generative models capture high-dimensional flow field distributions. Peebles et al.~\cite{peebles2023scalable} proposed the Diffusion Transformer (DiT) for global self-attention. Liu et al.~\cite{liu2022flow} proposed Rectified Flow for deterministic sampling via straight-line ODEs. Hu et al.~\cite{hu2019squeeze} proposed SE-Net to capture nonlinear couplings using channel attention.

\textbf{Improvements in this study:} Purely data-driven models may encounter non-physical local minima when subjected to strict conservation laws without proper constraints. In GeoFunFlow-3D (Sections 3.2 and 3.5), an approximate $2$-Wasserstein geodesic path is rebuilt in the latent space. Furthermore, we construct a $C^1$ continuous Variational Homotopy Scheduling and an exponential relaxation mechanism based on the numerical continuation method. This allows the model to smoothly shift from data topological alignment to strict Navier-Stokes conservation. It fundamentally establishes the global Well-posedness of the optimization problem.

\section{Methodology}

Before delving into the specific mathematical derivations, the global training pipeline is presented in \textbf{Algorithm \ref{alg:geofunflow}}. Let $\mathcal{P}$ denote the input 3D spatial point cloud, $\mathcal{F}_{in}$ represent the corresponding multi-modal geometric features, and $\theta$ parameterize the global generative neural network optimized over $E$ total epochs. During the flow matching evolution, $t \sim \mathcal{U}(0, 1)$ acts as the continuous time variable interpolating between the prior distribution and the target manifold. Furthermore, $\tau \in [0, 1]$ is defined as the normalized training progression step ($\tau = epoch/E$), which explicitly controls the dual optimization dynamics scheduling via the variational homotopy weight $\lambda(\tau)$.

\begin{algorithm}[htbp] 
\caption{GeoFunFlow-3D Training Pipeline} 
\label{alg:geofunflow}
\DontPrintSemicolon
\scriptsize  
\SetInd{0.3em}{0.3em} 
\SetAlgoLined

\KwIn{Point cloud $\mathcal{P}$, Features $\mathcal{F}_{in}$; Epochs $E$; Task domain $Task \in \{\Omega_{ext}, \Omega_{int}\}$}
\KwOut{Global network parameters $\theta^*$}
Initialize $\theta$\;
\For{$epoch = 1$ \KwTo $E$}{
    \tcc{1. Dual Optimization Dynamics Scheduling (Section 3.5)}
    $\tau \leftarrow epoch / E$\;
    $\lambda(\tau) \leftarrow \text{Cos\_Homo}(\tau, \tau_1, \tau_2)$ \tcp*{Variational Homotopy scheduling (Eq. 17)}
    $\lambda_{\mathrm{phys}}(\tau) \leftarrow \lambda_{\mathrm{base}} \cdot \eta^\tau$ \tcp*{Exponential physical relaxation (Eq. 18)}
    
    \For{each batch $(\mathcal{P}, \mathcal{F}_{in})$}{
        \tcc{2. Manifold Dimension Reduction and Geodesic Flow Matching (Sections 3.1 \& 3.2)}
        $Z_{\mathrm{geom}} \leftarrow \text{NW\_Kernel}(\text{GNO}(\mathcal{P}, \mathcal{F}_{in}))$ \tcp*{Density-invariant projection}
        $t \sim \mathcal{U}(0, 1)$, $Z_{\mathrm{noise}} \sim \mathcal{N}(0, I)$\;
        $Z_t \leftarrow (1 - t) Z_{\mathrm{noise}} + t Z_{\mathrm{geom}}$ \tcp*{W2 approximate geodesic line}
        $\hat{v}_\theta \leftarrow \text{DiT}(\text{SE}(Z_t), t)$ \tcp*{Direction prediction and feature calibration}
        $\mathcal{L}_{FM} \leftarrow \| \hat{v}_\theta - (Z_{\mathrm{geom}} - Z_{\mathrm{noise}}) \|_2^2$ \tcp*{Least action regression (Eq. 7)}
        
        \tcc{3. Frequency-domain Decoding and Real Gradient Pullback (Sections 3.2.3, 3.3 \& 3.4)}
        $u_{\mathrm{coarse}} \leftarrow \text{3D\_FNO}(\hat{v}_\theta)$ \tcp*{Spectral truncation continuous decoding}
        $u_{\mathrm{fine}} \leftarrow \mathcal{I}[u_{\mathrm{coarse}}] + M \cdot [\alpha \odot \mathcal{R}_{SATO}(\mathcal{Z}_{local})]$ \tcp*{SATO Fine Reconstruction (Eq. 14)}
        $\nabla_{\mathrm{real}}u_{\mathrm{coarse}} \leftarrow J_f^T \cdot \nabla_{\mathrm{lat}}^4(u_{\mathrm{coarse}})$ \tcp*{No-AD fourth-order real gradient (Eq. 10, 11)}
        $\mathcal{L}_{TV} \leftarrow \mathbb{E} \left[ \sum_{i \in \{x,y,z\}} \frac{|\nabla_i u_{\mathrm{coarse}}|}{R_i + \epsilon} \right]$ \tcp*{Anisotropic TV regularization (Eq. 9)}
        
        \tcc{4. Multi-modal Physical Hard Constraints (Sections 3.3.3 \& 3.4)}
        \uIf{ $Task \in \Omega_{ext}$ }{
            $M \leftarrow \sigma(\alpha \text{SDF} + \beta)$ \tcp*{Phase-field diffuse interface}
            $\vec{V}_{\mathrm{pred}} \leftarrow \text{Velocity}(u_{\mathrm{fine}})$ \tcp*{Extract predicted velocity field}
            $\vec{V}_{\mathrm{hard}} \leftarrow \vec{V}_{\mathrm{pred}} - \langle \vec{V}_{\mathrm{pred}}, \vec{n} \rangle \vec{n}$ \tcp*{Tangent bundle orthogonal projection (Eq. 15)}
            $\mathcal{R}_{PDE} \leftarrow \text{Navier\_Stokes}(\vec{V}_{\mathrm{hard}}, \nabla_{\mathrm{real}}u_{\mathrm{coarse}})$\;
            $\mathcal{L}_{\mathrm{phys}} \leftarrow \mathbb{E} [ M \cdot \|\mathcal{R}_{PDE}\|_2^2 \cdot |\det(J_f)| ]$ \tcp*{Singularity Cancellation (Eq. 13)}
        }
        \ElseIf{ $Task \in \Omega_{int}$ }{
            $M_{s} \leftarrow 1 - e^{-\gamma \|\nabla P\| / (P + \epsilon)}$ \tcp*{Phenomenological shock mask}
            $\mathcal{L}_{EOS} \leftarrow \mathbb{E} [ |P - \rho R T| ]$\;
            $\mathcal{L}_{th} \leftarrow \lambda_{iso} \mathbb{E}[\|\nabla S\|(1\!-\!M_{s})] + \lambda_{2nd} \mathbb{E}[\text{ReLU}(-\nabla P \!\cdot\! \nabla S) M_{s}]$\tcp*{Thermodynamic loss (Eq. 16)}
            $\mathcal{L}_{\mathrm{phys}} \leftarrow \mathcal{L}_{EOS} + \mathcal{L}_{th}$\;
        }
        
        \tcc{5. Dynamics Update}
        $\mathcal{L}_{\mathrm{total}} \leftarrow \mathcal{L}_{FM} + \mathcal{L}_{TV} + \lambda(\tau) \cdot \lambda_{\mathrm{phys}}(\tau) \cdot \mathcal{L}_{\mathrm{phys}}$ \tcp*{Action functional minimization (Section 4.5)}
        $\theta \leftarrow \theta - \eta \nabla_{\theta} \mathcal{L}_{\mathrm{total}}$ \tcp*{Variational path smooth sliding}
    }
}
\end{algorithm}

This section details the methodology of GeoFunFlow-3D, a physics-enhanced deep generative model designed to solve high-fidelity aerodynamic inference under complex 3D boundaries. The generation pipeline follows a strict computational logic. This includes spatial dimension reduction, geodesic evolution, discrete operators, multi-scale reconstruction, and dual scheduling.

The architecture contains five core stages. The mathematical proofs of convergence are in Section 4:

\begin{enumerate}
    \item \textbf{Geometry-Aware Encoding (Section 3.1):} Tangent bundle properties of point clouds are extracted. Physical fields are projected into a regular latent grid. This uses the partition of unity property of Nadaraya-Watson regression.
    \item \textbf{Latent Flow Matching (Section 3.2):} A flow matching path is built in the latent space. Physical dynamics are reconstructed via inverse Jacobian pullback and nonlinear feature calibration.
    \item \textbf{No-AD Discrete Engine (Section 3.3):} A central difference operator with $\mathcal{O}(h^4)$ error calculates physical residuals. This avoids AD-induced memory explosion and optimization oscillation.
    \item \textbf{Multi-scale Reconstruction (Section 3.4):} 3D-FNO performs continuous decoding. The SATO module compensates for residuals. Geometric phase-field and volumetric thermodynamic constraints are then applied.
    \item \textbf{Dual Optimization Dynamics (Section 3.5):} A homotopy scheduling protocol and an exponential relaxation mechanism are introduced to ensure smooth training evolution.
\end{enumerate}

\subsection{Geometry-Aware Encoding and Latent Grid Mapping}
Unstructured point clouds conflict with regular latent spaces. Thus, geometry-aware recognition is decoupled from operator dynamics evolution~\cite{wen2025geometry}. We construct a non-parametric geometric feature encoder.

\subsubsection{Multi-modal Geometric Embedding}
For an unstructured point cloud $\mathcal{P} = \{x_i\}_{i=1}^N \subset \mathbb{R}^3$, we construct a joint physical-geometric embedding vector $\mathcal{F}_{in}$. It extends Euclidean coordinates to a 9D space, including the unit normal vector $\vec{n}(x_i)$, signed distance field $\text{SDF}(x_i)$, and local principal curvature $\kappa(x_i)$:
\begin{equation}
\mathcal{F}_{in}(x_i) = \left[ x_i, \; \vec{n}(x_i), \; \text{SDF}(x_i), \; \kappa(x_i), \; \theta_{\mathrm{dummy}} \right]^T \in \mathbb{R}^9 .
\end{equation}
This embedding effectively endows the discrete points with the tangent bundle properties of a Riemannian manifold.

\subsubsection{Neural Integral Operator Approximation}
The point cloud is modeled as a $k$-nearest neighbor graph $\mathcal{G} = (\mathcal{V}, \mathcal{E})$. Based on graph neural operator (GNO) theory, it is discretized into a message passing system:
\begin{equation}
h^{(l+1)}(x_i) = \sigma \left( W h^{(l)}(x_i) + \frac{1}{|\mathcal{N}(x_i)|} \sum_{x_j \in \mathcal{N}(x_i)} \mathcal{K}_\theta(x_i - x_j) \odot h^{(l)}(x_j) \right) .
\end{equation}
The parameterized kernel $\mathcal{K}_\theta$ extracts spatial topological features while maintaining translation invariance.

\subsubsection{Density-Invariant Projection Operator}
To map the unordered set $\mathcal{H}_{point}$ to the regular latent grid $\Xi_{\mathrm{lat}}$, we introduce Nadaraya-Watson (NW) kernel regression (see Algorithm \ref{alg:geofunflow}, line 7) and derive an inverse distance projection mechanism:
\begin{equation}
\tilde{\phi}(\xi_k) = \text{LN} \left( \sum_{x_j \in \mathcal{N}(\xi_k)} \frac{\omega(d(\xi_k, x_j))}{\sum_{x_m \in \mathcal{N}} \omega(d(\xi_k, x_m))} \cdot \left[ \mathcal{K}_\phi(\xi_k - x_j) \odot h^{(L)}(x_j) \right] \right) .
\end{equation}
Here, $\xi_k$ is the grid anchor, and $\omega(d) = 1/(d + \epsilon)$ is the kernel weight. This constructs a "partition of unity" via weight normalization. It ensures robustness to point cloud density. The Lipschitz continuity proof is in \ref{App:Lipschitz}.

\subsubsection{Feature Auto-Encoder (FAE) Warm-up}
The latent grid $\Xi_{\mathrm{lat}}$ must restore the original physical information. Thus, we introduce an FAE warm-up strategy via an end-to-end reconstruction task.

Assume the target physical field (e.g., pressure $C_p$) on $\mathcal{P}$ is $\mathcal{U} = \{u(x_i)\}_{i=1}^N$. The FAE encoding mapping $\mathcal{E}_\theta$ combines GNO and NW projection:
\begin{equation}
\tilde{\phi} = \mathcal{E}_\theta(\mathcal{P}, \mathcal{U}) \in \mathbb{R}^{D_{\mathrm{lat}} \times H \times W \times D} .
\end{equation}
A continuous decoding operator $\mathcal{D}_\psi$ (3D-FNO, Section 3.4.1) maps features back to the physical space to get $\hat{\mathcal{U}} = \{\hat{u}(x_i)\}_{i=1}^N$:
\begin{equation}
\hat{u}(x_i) = \mathcal{D}_\psi(\tilde{\phi}, x_i) .
\end{equation}

During warm-up, parameters $\{\theta, \psi\}$ are optimized using the reconstruction loss $\mathcal{L}_{FAE}$:
\begin{equation}
\mathcal{L}_{FAE} = \frac{1}{N} \sum_{i=1}^N \|u(x_i) - \hat{u}(x_i)\|_2^2 + \lambda_{reg} \|\tilde{\phi}\|_2^2 .
\end{equation}

\begin{figure}[htbp]
    \centering
    \includegraphics[width=0.8\textwidth]{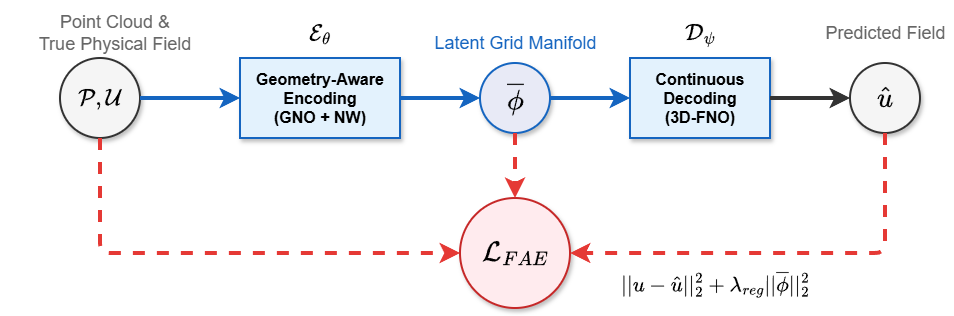}
    \caption{The computational graph of the Feature Auto-Encoder (FAE) warm-up stage. Solid arrows indicate the forward data mapping from the physical domain to the latent grid and back to the predicted field. Dashed red arrows represent the backward pathway for loss aggregation and parameter optimization.}
    \label{fig:fae_warmup}
\end{figure}

To provide an intuitive understanding of this process, Figure \ref{fig:fae_warmup} illustrates the computational graph of the FAE warm-up stage. The solid blue and black arrows delineate the forward data-driven pathway. The initial point cloud and true physical field ($\mathcal{P}, \mathcal{U}$) are processed through the geometry-aware encoder ($\mathcal{E}_\theta$) and projected onto the latent grid manifold ($\tilde{\phi}$). Subsequently, the continuous decoder ($\mathcal{D}_\psi$) reconstructs the predicted field ($\hat{u}$). Conversely, the dashed red arrows represent the backward loss aggregation and optimization pathway. The $\mathcal{L}_{FAE}$ node explicitly gathers the end-to-end prediction error while extracting information from the latent grid to apply the regularization penalty.

The second term in the loss function prevents latent variable scattering. This locks the homeomorphic mapping parameters. It provides a dense, physically interpretable initial manifold for 2-Wasserstein ($W_2$) flow matching (Section 3.2). Notably, this FAE warm-up phase learns the mapping from geometry to the latent space without relying on any downstream physical boundary conditions. As a result, the forward generation phase remains a strict zero-observation process.

\subsection{Generative Flow Matching and Dynamics Reconstruction}
Fitting the evolution operator in a high-dimensional Euclidean space triggers the "curse of dimensionality". Inspired by manifold adaptivity~\cite{kumar2026flow}, we build a flow matching mechanism in the low-dimensional latent space $\Xi_{\mathrm{lat}} \subset \mathbb{R}^d$. Mathematical proofs of its error bound and geodesic equivalence are detailed in Section 4.2.

\subsubsection{Construction of Conditional Flow Matching Objective}
In $\Xi_{\mathrm{lat}}$, GeoFunFlow-3D learns a time-varying vector field $v_t$ (see Algorithm \ref{alg:geofunflow}, line 10). This field pushes the prior Gaussian noise $p_0$ to the target flow field $p_1$. We choose a Gaussian linear interpolation as the conditional probability path ($z_t = (1-t) z_0 + t z_1$). Thus, the flow matching regression loss is defined as:
\begin{equation}
\mathcal{L}_{FM}(\theta) = \mathbb{E}_{t, z_0, z_1} \left\| \hat{v}_\theta(z_t, t) - (z_1 - z_0) \right\|_2^2 .
\end{equation}

\subsubsection{Subspace Projection in the Tangent Vector Space}
We introduce the physical constraint $\mathcal{F}(u)=0$. Mathematically, we perform a physical subspace projection on the tangent vector $v_t$ at each time step $t$. This forces the predicted direction $\hat{v}_\theta$ to stay within a specific tangent bundle subset. This subset satisfies the Navier-Stokes solution manifold.

\subsubsection{Inverse Jacobian Physical Reconstruction and Anisotropic Total Variation Regularization}
The latent evolution must be mapped back to the 3D physical space $\Omega_{\mathrm{real}}$. We define an anisotropic affine homeomorphic mapping $\xi = f(X)$ and reconstruct the real physical velocity via an inverse Jacobian ($J_f^{-1}$) pullback:
\begin{equation}
V_{\mathrm{real}}(X) = J_f^{-1} v_t(\xi) = \text{diag}(R_x+\epsilon, R_y+\epsilon, R_z+\epsilon) \cdot (z_1 - z_0) .
\end{equation}
Furthermore, numerical oscillations in transonic shock regions must be suppressed. An anisotropic 3D total variation regularization functional (TV Loss) is derived based on this pullback:
\begin{equation}
\mathcal{L}_{TV} = \mathbb{E} \left[ \frac{|\nabla_x \tilde{\phi}|}{R_x} + \frac{|\nabla_y \tilde{\phi}|}{R_y} + \frac{|\nabla_z \tilde{\phi}|}{R_z} \right] .
\end{equation}
The spatial bounding radius $R$ acts as a Jacobian feature multiplier. It serves as an auxiliary flux limiter.

\subsubsection{Nonlinear Feature Adaptive Calibration}
We embed a feature adaptive calibration module into the core path. It captures the nonlinear coupling among density, velocity, and pressure. First, we extract a macroscopic physical state descriptor $z_c = \frac{1}{|\Xi_{\mathrm{lat}}|} \int_{\Xi_{\mathrm{lat}}} u_c(\xi) d\xi$ on the manifold. Then, we derive a gated dependency matrix $s = \sigma(W_2 \delta(W_1 z))$ via an information bottleneck. Finally, we perform a dynamic recalibration $\tilde{U}_c = s_c \cdot u_c$. This adaptively amplifies high-frequency features like shock gradients.

\subsection{Jacobian-Aware Discrete Engine and Residual Assembly}

To resolve the severe memory explosion and Neural Tangent Kernel (NTK) spectral bias caused by automatic differentiation (AD) when computing high-order derivatives, we introduce a high-order discrete difference engine in the latent space (see Algorithm \ref{alg:geofunflow}, lines 14).

\subsubsection{Tangent Bundle Pullback from Physical to Latent Space}
To enable stable evolution on regular grids, we construct an affine homeomorphic transformation $\xi = f(X) = (X - c) \oslash (R + \epsilon)$. According to the chain rule, we pull back the real gradient $\nabla_{\mathrm{real}}$ from the latent gradient $\nabla_{\mathrm{lat}}$ via the Jacobian transpose:
\begin{equation}
\nabla_{\mathrm{real}}\phi(X) = J_f^T \nabla_{\mathrm{lat}}\tilde{\phi}(\xi) \implies \frac{\partial \phi}{\partial x_i} = \frac{1}{R_i + \epsilon} \frac{\partial \tilde{\phi}}{\partial \xi_i}, \quad i \in \{x, y, z\} .
\label{eq:jacobian_pullback}
\end{equation}
This decoupled operation uses a diagonal inverse Jacobian multiplier. It lays the foundation for No-AD differentiation. However, 3D geometries often exhibit extreme curvature mutations or severe distortions. The inverse Jacobian matrix $J_f^{-1}$ is highly susceptible to singularity ($\det(J_f) \to 0$). This causes non-physical high-frequency oscillations in the pulled-back gradients. To suppress this issue algebraically, the parameter $\epsilon$ in Equation \eqref{eq:jacobian_pullback} acts as a diagonal Tikhonov Regularization (see \ref{app:tikhonov} for detailed sensitivity analysis and hyperparameter selection). It constructs a generalized pseudo-inverse by truncating tiny singular values. This ensures the gradient trajectory remains in a Lipschitz continuous stable region even under extreme geometries~\cite{hansen1987rank}. A Krylov subspace preconditioning mechanism can be introduced later for non-diagonal singular mappings.

\subsubsection{Engineering Implementation of High-Fidelity Discrete Differential Operators}
The operator $\nabla_{lat}$ must be explicitly calculated on the latent grid. Traditional AD can introduce spectral distortion under certain high-frequency conditions. (The rigorous proof demonstrating that this discrete formulation bounds the NTK spectral explosion is detailed in Theorem 1, Section 4.3). Therefore, we opt for a high-fidelity fourth-order central difference operator derived via Taylor's expansion. Its truncation error is $\mathcal{O}(h^4)$ (derivations are detailed in \ref{App:Taylor}):
\begin{equation}
\tilde{\phi}'(\xi) = \frac{-\tilde{\phi}(\xi + 2h) + 8\tilde{\phi}(\xi + h) - 8\tilde{\phi}(\xi - h) + \tilde{\phi}(\xi - 2h)}{12h} .
\end{equation}
However, some nodes are close to boundaries. Complete $\pm 2h$ symmetric neighboring nodes are unavailable. Applying central differences here causes out-of-bounds errors. To avoid boundary accuracy degradation, ghost cells are excluded, which prevents non-physical extrapolation artifacts. Instead, we introduce a high-order one-sided difference mechanism to maintain the $\mathcal{O}(h^4)$ error. For example, we apply a forward five-point scheme at the left boundary node $\xi_{bound}$:
\begin{equation}
\begin{aligned}
\tilde{\phi}'(\xi_{bound}) &= \frac{1}{12h} \Big[ -25\tilde{\phi}(\xi) + 48\tilde{\phi}(\xi + h)- 36\tilde{\phi}(\xi + 2h) + 16\tilde{\phi}(\xi + 3h) - 3\tilde{\phi}(\xi + 4h) \Big] .
\end{aligned}
\end{equation}
A symmetric backward scheme is applied at the right boundary. Difference templates are dynamically switched at the tensor edges. 

\begin{figure}[htp]
    \centering
    \includegraphics[width=0.8\textwidth]{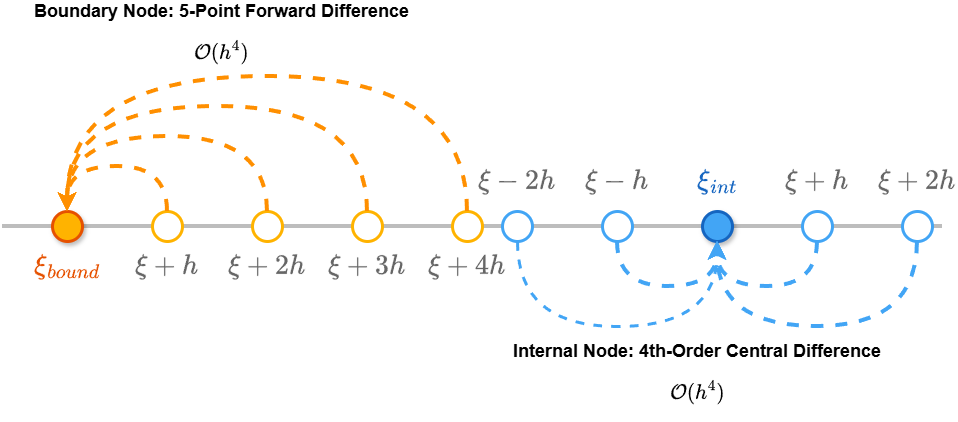}
    \caption{Illustration of the No-AD discrete differential mechanism. The blue dashed arrows represent the symmetric 4th-order central difference for internal nodes. The orange dashed arrows illustrate the 5-point forward difference for boundary nodes, dynamically avoiding out-of-bounds errors while globally maintaining a consistent $\mathcal{O}(h^4)$ accuracy.}
    \label{fig:no_ad_stencil}
\end{figure}

To visually clarify this dynamic stencil switching, Figure \ref{fig:no_ad_stencil} illustrates our discrete differential mechanism. The blue components represent an internal node ($\xi_{int}$). The dashed blue arrows show the standard fourth-order central difference, which gathers information from symmetric neighboring nodes ($\pm h, \pm 2h$). Conversely, the orange components highlight a boundary node ($\xi_{bound}$). To prevent out-of-bounds errors, the dashed orange arrows depict the five-point forward difference scheme. This one-sided scheme extracts information strictly from the boundary itself and the adjacent interior nodes (up to $\xi+4h$). Both pathways explicitly maintain an $\mathcal{O}(h^4)$ truncation error. This visually confirms that the framework ensures consistent global accuracy without needing artificial ghost cells.

Thus, GeoFunFlow-3D constructs a closed computational graph algebraically. This maintains a consistent $\mathcal{O}(h^4)$ error globally and effectively mitigates boundary gradient divergence. While high-order one-sided schemes can inherently introduce numerical instability, our framework mitigates this risk through its formulation. Specifically, the variational functional (Eq. 13) and the TV loss (Eq. 9) act as implicit numerical filters. They effectively absorb and suppress high-frequency oscillations originating from the boundaries, thereby maintaining stable gradient updates and facilitating robust global convergence.

\subsubsection{Variational Physical Functional Assembly and Singularity Cancellation of Discrete Residuals}
Classical PINNs rely on strong-form PDE residuals based on pointwise Jacobian mappings. This requires high function smoothness and may face convergence challenges near geometric singularities. To mitigate potential numerical instability from strong local gradients, we transition from the pure pointwise strong form to a Lipschitz continuous variational functional integral:
\begin{equation}
\mathcal{L}_{\mathrm{phys}} = \mathbb{E}_{\Xi_{\mathrm{lat}}} \left[ M(\xi) \cdot \left\| \mathcal{F}(J_f^T \nabla_{\mathrm{lat}}\tilde{\phi}, \tilde{\phi}) \right\|_p^p \right] .
\end{equation}
A Hadamard product is executed with the implicit Lagrangian phase-field mask $M(\xi)$. This zeros out invalid residuals in non-fluid domains. Mathematically, this expected integral expression naturally resists distortion. At geometric singularities, the inverse matrix $J_f^{-T}$ approaches infinity. Meanwhile, the Jacobian determinant $|\det(J_f)|$ in the integral measure approaches zero. These two factors undergo algebraic cancellation and physical smoothing within the integral kernel~\cite{kharazmi2019vpinns, ryck2022wpinns}. This fundamentally eliminates the impact of grid distortion on the Hessian matrix. It achieves robustness similar to weak-form PINNs (wPINNs).

\subsection{Multi-scale Reconstruction and Physical Hard-Constraining}

\subsubsection{Global Frequency-Domain Continuous Decoding Based on 3D-FNO}
Based on the Fourier Neural Operator, latent space feature evolution is mapped as: $(\mathcal{K}_{\phi} v)(x) = \mathcal{F}^{-1} \left( R_{\phi} \cdot (\mathcal{F} v) \right)(x)$. We perform a strict spectral truncation in the frequency domain, actively filtering high-frequency components. This promotes global topological smoothness.

\subsubsection{Residual Compensation Ansatz of the SATO Module}
\begin{figure}[htbp] 
    \centering 
    \includegraphics[width=1.0\linewidth]{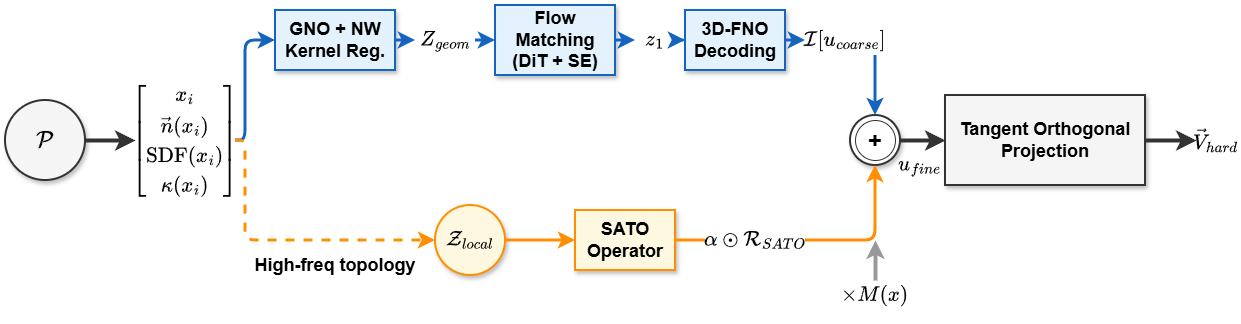}
    \caption{Mathematical computational graph of the SATO module's multi-scale residual compensation process. \textbf{Color-coded arrow semantics}: Solid blue arrows delineate the macroscopic flow topology pathway via 3D-FNO decoding and global interpolation. Orange arrows (dashed for extraction, solid for propagation) represent the microscopic high-frequency residual compensation pathway generating basis functions ($\mathcal{R}_{SATO}$). The grey arrow indicates the application of the phase-field mask $M(x)$. The pathways are merged using an addition operator ($\bigoplus$) to produce the fine-scale field $u_{\mathrm{fine}}$, followed by a tangent orthogonal projection (dark grey arrows) to ensure local physical consistency near sharp gradients such as shock waves.} 
    \label{fig:sato} 
\end{figure}

Frequency filtering intrinsically weakens local physical extremes in boundary layers. Thus, we introduce a topology-aware super-resolution module (SATO), establishing a strong connection between macroscopic flow topologies and microscopic physical features. As rigorously depicted in the mathematical computational graph (Figure \ref{fig:sato}), the forward inference pipeline splits into two synergistic pathways to reconstruct the flow field from the multi-modal geometric embedding.

First, the upper branch (indicated by \textbf{solid blue arrows}) utilizes the 3D-FNO output to provide a macroscopic smooth background ($u_{\mathrm{coarse}}$). This smooth field undergoes a global interpolation to project it onto the high-resolution target grid, represented by the term $\mathcal{I}[u_{\mathrm{coarse}}](x)$. 

Simultaneously, to recover the lost high-frequency details, the lower branch routes fine-scale geometric topologies ($\mathcal{Z}_{local}$) via a \textbf{dashed orange arrow} into the SATO Operator. It extracts fine-scale geometric features to generate the residual compensation field (propagated via a \textbf{solid orange arrow}), mathematically serving as neural microscopic basis functions $\mathcal{R}_{SATO}$. These two pathways converge at an addition operator ($\bigoplus$) to produce the fine-scale reconstruction. Crucially, the SATO residual is multiplied by the implicit phase-field mask $M(x)$, indicated by the \textbf{grey arrow}. This mechanism effectively bypasses the resolution limitations typically imposed by the spatial step $h$, clearly capturing sharp physical changes like transonic shocks without requiring massive latent grid refinement. The final evolution ansatz for the fine-scale flow field is defined as:
\begin{equation}
u_{\mathrm{fine}}(x) = \mathcal{I}[u_{\mathrm{coarse}}](x) + M(x) \cdot \left[ \alpha \odot \mathcal{R}_{SATO}(x; \mathcal{Z}_{local}) \right] .
\end{equation}
Subsequently, $u_{\mathrm{fine}}$ undergoes a tangent orthogonal projection (indicated by \textbf{dark grey arrows}) to strictly enforce the Euler wall tangency condition ($\vec{V}_{\mathrm{hard}}$). This module acts as a neural generalization of the multiscale finite element method (MsFEM) (the mathematical proof regarding this spatial local consistency is detailed in Section 4.4).

\subsubsection{Surface Hard-Constraining Based on Phase-Field Theory and Tangent Bundle Projection}
We apply two boundary hard constraints:

\textbf{1. Implicit Lagrangian Phase-Field Mask:} We construct a diffuse interface mask $M(x) = \sigma(\alpha \cdot \mathcal{D}(x) + \beta)$. It locks the non-penetration constraint within the boundary layer. Defective CAD grids (e.g., local holes) pose challenges for the robustness of standard SDF masks. Thus, we integrate a topology repair mechanism. It uses level-set methods or Gaussian smoothing to fix geometric singularities. This prevents hard constraint failures.

\textbf{2. Tangent Space Orthogonal Projection:} We construct a linear projection operator and subtract the normal velocity component to enforce the Euler wall condition:
\begin{equation}
\vec{V}_{\mathrm{hard}} = \vec{V}_{\mathrm{pred}} - \langle \vec{V}_{\mathrm{pred}}, \vec{n} \rangle \vec{n} .
\end{equation}

\subsubsection{Multi-modal Decoupling and Non-isentropic Volumetric Thermodynamic Constraints}
For highly compressible internal flows, the generated fields must satisfy the ideal gas equation ($\mathcal{L}_{EOS} = \mathbb{E}_{\Omega} [|P - \rho R T|]$). We construct a phenomenological shock mask $M_{\mathrm{shock}} = 1 - \exp( -\gamma ||\nabla P|| / (P + \epsilon) )$. Then, we introduce a local second law of thermodynamics functional:
\begin{equation}
\begin{split}
\mathcal{L}_{thermo} &= \lambda_{iso} \mathbb{E} \left[ ||\nabla S|| \cdot (1 - M_{\mathrm{shock}}) \right] \\
&\quad + \lambda_{2nd} \mathbb{E} \left[ \text{ReLU}(-\nabla P \cdot \nabla S) \cdot M_{\mathrm{shock}} \right] .
\end{split}
\end{equation}
In shock regions ($M_{\mathrm{shock}} \to 1$), the ReLU function strictly prohibits local entropy decrease. This phenomenological penalty is critical for complex transonic internal flow fields (such as the Rotor37 compressor); it strictly enforces the Second Law of Thermodynamics, facilitating the robust capture of 3D detached shock topologies while effectively suppressing numerical oscillations.

\subsection{Implementation of Dual Optimization Dynamics}
We design a dual scheduling system to balance physical constraints and data manifolds. (The theoretical proof of how this scheduling ensures path smoothing and global Well-posedness is detailed in Section 4.5.)

First, a Lagrangian multiplier $\lambda(\tau)$ is defined. It acts as a $C^1$ continuous cosine annealing homotopy function for the training step $\tau$:
\begin{equation}
\lambda(\tau) = 
\begin{cases} 
0, & \tau < \tau_1 \\
\frac{\lambda_{max}}{2} \left[ 1 - \cos\left( \frac{\pi(\tau - \tau_1)}{\tau_2 - \tau_1} \right) \right], & \tau_1 \le \tau \le \tau_2 \\
\lambda_{max}, & \tau > \tau_2
\end{cases} .
\end{equation}

\begin{figure}[htbp] 
    \centering 
    \includegraphics[width=0.7\linewidth]{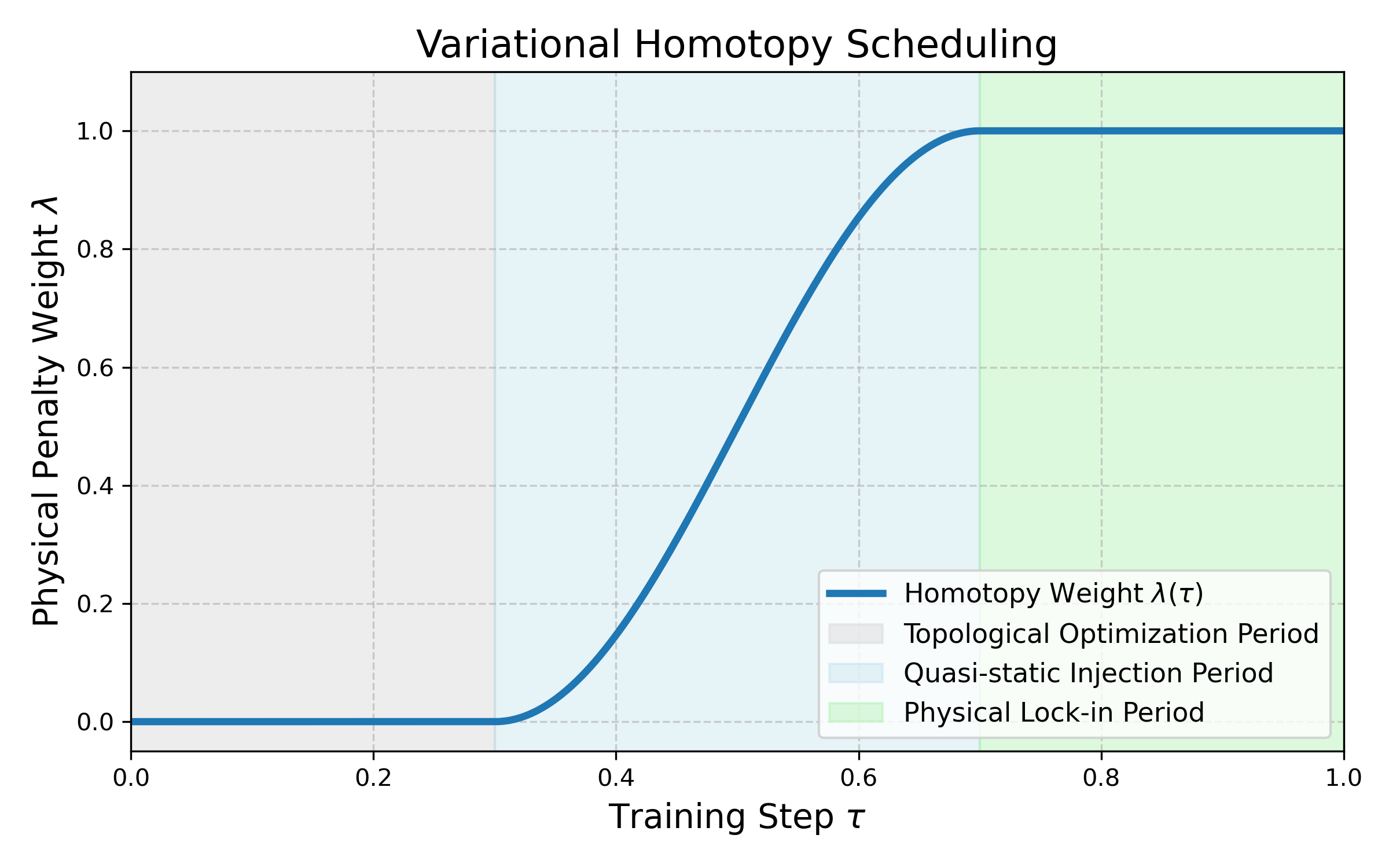} 
    \caption{The evolution curve of the Variational Homotopy Scheduling weight $\lambda(\tau)$ across training steps, visually delineating the topological optimization period, quasi-static injection period, and physical lock-in period.} 
    \label{fig:homotopy} 
\end{figure}

The system experiences three stages: the topological optimization period ($\lambda=0$), the quasi-static injection period, and the steady-state convergence lock ($\lambda=\lambda_{max}$), as strategically delineated in Figure \ref{fig:homotopy}. This protocol effectively partitions the training process into three distinct functional regimes to navigate the non-convex optimization landscape of high-dimensional fluid dynamics. During the initial ``Topological Optimization Period'' ($\tau < \tau_1$), the penalty weight remains at zero, allowing the generative flow matching model to prioritize the alignment of macroscopic geometric topologies and data manifolds without the interference of rigid physical constraints. As training progresses into the ``Quasi-static Injection Period'' ($\tau_1 \le \tau \le \tau_2$), the weight $\lambda$ follows a $C^1$ continuous cosine annealing path. This smooth transition acts as a numerical continuation method, adiabatically injecting the non-convex Navier-Stokes constraints into the optimization objective to prevent sudden gradient conflicts and optimizer instability. Finally, in the ``Physical Lock-in Period'' ($\tau > \tau_2$), $\lambda$ is sustained at its maximum value $\lambda_{max}$. During this stage, the model is strictly constrained by the volumetric thermodynamic and momentum residuals, ensuring the precise calibration of high-frequency physical features, such as detached shocks and boundary layers, while maintaining global convergence.

Second, we apply an exponential physical relaxation protocol to the ODE trajectory generation. This releases network capacity to fit high-frequency data mutations later in the evolution:
\begin{equation}
\lambda_{\mathrm{phys}}(\tau) = \lambda_{\mathrm{base}} \cdot \eta^\tau \quad (\eta \in (0,1)) .
\end{equation}

\section{Theoretical Analysis of Convergence}

\subsection{Theoretical Summary: Trinity Convergence Guarantee of GeoFunFlow-3D}

\begin{figure}[htbp] 
    \centering 
    \includegraphics[width=0.9\linewidth]{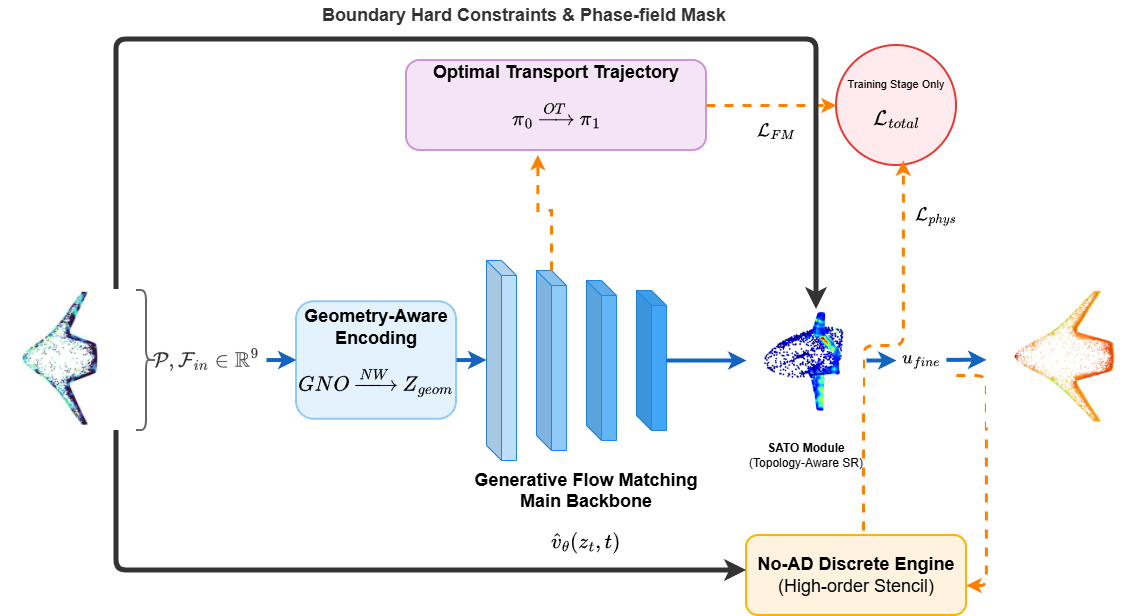} 
    \caption{The mathematical computational graph and trinity collaborative architecture of GeoFunFlow-3D. \textbf{Color-coded arrow semantics}: Solid dark blue arrows denote the forward data flow for zero-observation inference; red dashed arrows represent the loss-driven optimization flow for training; and orange dashed arrows illustrate the coupling between the generated fine-scale field ($u_{fine}$) and the No-AD discrete engine for physical residual calculation. The top grey arc indicates the enforcement of boundary hard constraints and phase-field masks across the spectral and spatial modules.} 
    \label{fig:trinity} 
\end{figure}

The architecture of GeoFunFlow-3D is constructed upon a rigorous mathematical computational graph that coordinates temporal, spectral, and spatial dimensions to ensure physical consistency. As illustrated in Figure~\ref{fig:trinity}, the framework employs a semantic color-coding scheme to distinguish between different operational modes. The \textbf{solid dark blue arrows} delineate the forward inference pipeline, representing the deterministic mapping from multi-modal geometric embeddings to high-fidelity aerodynamic fields. Conversely, the \textbf{red dashed arrows} indicate the backward loss aggregation and optimization pathways, which are active exclusively during the training stage for minimizing the total action functional. Furthermore, the \textbf{grey arc arrows} represent the global intervention of boundary hard constraints and phase-field masks, ensuring that the generated flow fields adhere strictly to non-penetration conditions.

We sequentially conduct rigorous theoretical analyses below, covering error bounds, linear spectral stability, multi-scale mapping equivalence, and dynamic global Well-posedness.

\subsection{Latent Space Homeomorphic Mapping and Optimal Transport Geodesic Energy Bound}

Fluid mechanics evolution is strictly controlled by a few macroscopic parameters (e.g., Mach number, angle of attack). Thus, the real physical distribution $\mathbb{P}_{\mathcal{X}}$ is supported on a low-dimensional smooth manifold $\mathcal{M}_{\mathrm{phys}}$. Its intrinsic dimension is $d$ ($d \ll D$). According to the Minimax Lower Bound principle, the dimension reduction strictly compresses the irreducible approximation error for the velocity field. The error is reduced from $\mathcal{O}(N^{-1/D})$ to $\mathcal{O}(N^{-1/d})$.

Furthermore, by solving the ordinary differential evolution equation and introducing Grönwall's Inequality, we mathematically prove the global generation error bound:
\begin{equation}
W_2(\mathbb{P}_{gen}, \mathbb{P}_{data}) \le \mathcal{O}\left( N^{-\frac{1}{d}} \right) \cdot e^L .
\end{equation}
The exponential term $e^L$ reveals the time-integrated explosion of tiny errors during inference. The error base is controlled via dimension reduction. Thus, this framework mathematically mitigates the exponential disaster of differential propagation (derivations are in \ref{App:Gronwall}).

In addition, flow matching efficiently approximates dynamic optimal transport (Dynamic OT). The path transferring distribution $p_0$ to $p_1$ with minimum physical kinetic energy is obtained by minimizing the Benamou-Brenier functional~\cite{benamou2000computational}:
\begin{equation}
\min_{v_t, p_t} \int_0^1 \int_{\Xi_{\mathrm{lat}}} \frac{1}{2} \|v_t(z)\|^2 p_t(z) dz dt, \quad \text{s.t. } \frac{\partial p_t}{\partial t} + \nabla \cdot (p_t v_t) = 0 .
\end{equation}
The constraint is the mass conservation equation. According to the Euler-Lagrange Equation, minimum kinetic energy requires zero material acceleration ($dv_t / dt = 0$). Thus, Gaussian linear interpolation is selected. Its trajectory acts as an approximate geodesic under the $2$-Wasserstein distance ($W_2$)~\cite{lipman2023flow}. This provides a theoretically optimal transport path with minimized kinetic energy.

\subsection{NTK Spectral Boundedness Analysis of the No-AD Discrete Differential Engine}

Linear stability is established based on the No-AD operator. Traditional PINNs rely on automatic differentiation (AD) for high-order derivatives. Assume the network contains a $k$-th order spatial derivative. Under AD, the Fourier symbol of the continuous operator is $(i\omega)^k$. According to Neural Tangent Kernel (NTK) theory~\cite{wang2022when}, fitting a frequency $\omega$ amplifies the corresponding eigenvalue $\lambda_{AD}(\omega)$:
\begin{equation}
\lambda_{AD}(\omega) \propto |(i\omega)^k|^2 = \omega^{2k} .
\end{equation}
For high-frequency details ($\omega \to \infty$), the maximum NTK eigenvalue explodes unboundedly at $\mathcal{O}(\omega^{2k})$. This can lead to severe Hessian Rigidity and ill-conditioned oscillations.

Comparing continuous and discrete frequency-domain symbols yields the core theorem:

\begin{theorem}[Discrete NTK Spectral Boundedness]
Unlike unbounded AD eigenvalues, the maximum NTK eigenvalue of the No-AD fourth-order discrete operator is strictly bounded by the spatial step $h$.
\end{theorem}

\begin{proof}
Applying the fourth-order difference operator to the Fourier basis $e^{i\omega \xi}$ yields the mapping symbol $\hat{D}_4(\omega) = i \left( \frac{8\sin(\omega h) - \sin(2\omega h)}{6h} \right)$. Since $|\sin(x)| \le 1$, a strict upper bound exists for any $\omega$:
\begin{equation}
|\hat{D}_4(\omega)| \le \frac{8|\sin(\omega h)| + |\sin(2\omega h)|}{6h} \le \frac{9}{6h} = \frac{3}{2h} .
\end{equation}
Therefore, the eigenvalue growth for the No-AD matrix is algebraically bounded:
\begin{equation}
\lambda_{No-AD}(\omega) \propto |\hat{D}_4(\omega)|^{2k} \le \left( \frac{3}{2h} \right)^{2k} = \mathcal{O}(h^{-2k}) .
\end{equation}
\end{proof}

The theoretical advantage of the No-AD engine is visually evidenced in Figure \ref{fig:spectrum}, which compares the spectral response of the continuous AD operator against our discrete implementation. In the log-scaled eigenvalue magnitude plot, the traditional AD operator (represented by the red dashed line) exhibits a monotonic and unbounded explosion in eigenvalue magnitude as the frequency $\omega$ increases. This phenomenon leads to severe stiffness in the optimization landscape, where high-frequency errors (noise) dominate the gradient descent process and trigger numerical instability. In stark contrast, the No-AD discrete engine (represented by the green solid line) exhibits an inherent band-limited filtering property. While it accurately tracks the physical signal in the low-frequency regime, it forcibly truncates the response in high-frequency regions. This ensures that the maximum eigenvalue is strictly upper-bounded by the spatial grid resolution $h$, as proven in Theorem 1. 

\begin{figure}[htbp] 
    \centering 
    \includegraphics[width=0.75\linewidth]{fig_spectrum.png} 
    \caption{Comparison of the spectral response between continuous automatic differentiation (AD) and the No-AD discrete operator. The No-AD mechanism acts as an ideal band-limited filter, effectively truncating the eigenvalue explosion ($\mathcal{O}(\omega^{2k})$) to eliminate high-frequency spectral bias and optimization oscillations.} 
    \label{fig:spectrum} 
\end{figure}

This reveals the dynamic essence of the No-AD operator. It acts as an ideal band-limited filter. High-frequency responses are smoothly truncated to $\mathcal{O}(1/h)$. This prevents the $\omega^{2k}$ explosion and establishes strict numerical Well-posedness for high-fidelity 3D aerodynamic inference.

\subsection{Spatial Local Consistency Proof of the SATO Module and MsFEM Mapping}

Frequency-domain operators often exhibit reduced accuracy near strong boundary layer gradients. Thus, SATO is introduced to enforce local physical consistency. Rather than functioning as a standard spatial interpolator, SATO serves as a nonlinear generalization of the multiscale finite element method (MsFEM) in the neural function space.

Classical MsFEM constructs microscopic basis functions $\phi_i$ on a coarse grid. These satisfy the local elliptic problem $\mathcal{L} \phi_i = 0$~\cite{hou1997multiscale}. In SATO, the macroscopic smooth field represents the coarse grid solution. "Neural microscopic basis functions" $\mathcal{R}_{SATO}$ are dynamically generated from high-frequency topological features. Its optimization strictly minimizes the local physical residual:
\begin{equation}
\min \| \mathcal{L} (\mathcal{I}[u_{\mathrm{coarse}}] + \mathcal{R}_{SATO}) \|_{\Omega_{local}}^2 .
\end{equation}
This is isomorphic to the MsFEM variational principle. Crucially, SATO automatically adjusts activation weights in high-gradient regions. This corresponds to the MsFEM over-sampling technique. It effectively eliminates resonance errors between coarse and fine grids.

\subsection{Global Convergence Analysis Driven by Variational Homotopy Scheduling}

Fluid mechanics equations are non-convex and contain local minima traps. Thus, we introduce Variational Path Smoothing and formulate flow generation as a PDE-constrained functional minimization. The generalized action functional is $\mathcal{S}(u; \tau) = \mathcal{J}_{data}(u) + \lambda(\tau) \cdot \| \mathcal{F}(u) \|_2^2$.

Under static strong coupling (constant $\lambda$), the condition number of the Hessian matrix $\mathbf{H} = \mathbf{H}_{data} + \lambda \mathbf{H}_{\mathrm{phys}}$ approaches infinity ($\kappa(\mathbf{H}) \to \infty$)~\cite{wang2021understanding}. This can introduce severe optimizer oscillations in high-curvature directions.

To resolve this Hessian Rigidity, we construct a continuous mapping $H(u, \gamma)$ using the homotopy function $\lambda(\tau)$. This follows the Numerical Continuation Method~\cite{allgower2012numerical}. It induces the optimization path to slip adiabatically within the adjoint space. The $C^1$ continuous cosine function guarantees smooth injection without gradient shocks. This ensures a smooth evolution from pure data-driven topological alignment to the Navier-Stokes solution space. Ultimately, it establishes the global Well-posedness of the nonlinear optimization problem.

\section{Experiments}
Comprehensive evaluations were conducted on two industrial-grade datasets. These are the BlendedNet dataset for aerodynamic surface design and the NASA Rotor37 dataset for 3D transonic thermodynamic volumetric fields.

\subsection{Experimental Setup and Baselines}

The BlendedNet dataset~\cite{sung2025blendednet} contains 8,830 aerodynamic surface samples of blended-wing-body configurations. The ground truth data were generated by solving the Reynolds-Averaged Navier-Stokes (RANS) equations closed with the Spalart-Allmaras turbulence model. Each case employed a high-fidelity mesh consisting of 9 to 14 million volume cells with $y^+ < 1$ to resolve the boundary layer. Flight conditions were determined via Latin Hypercube Sampling (LHS), covering a Mach number range of $0.05 \sim 0.5$, an angle of attack range of $-10^\circ \sim 20^\circ$, and a Reynolds length range of $0.1 \sim 10$ m. We adopted mean absolute error (MAE) and relative $L_2$ error as metrics, and set up a few-shot test environment with only 100 and 500 samples. This verified the generalization capability under extreme data scarcity.

The Rotor37 dataset contains 1,200 dense 3D volumetric field samples inside a transonic compressor rotor, sourced from the open-access PLAID 3D RANS database. The simulations were conducted at the design rotational speed of 17,188.7 RPM with standard sea-level inlet total pressure (101,325 Pa) and total temperature (288.15 K). Each single-passage sample utilized a high-resolution mesh of approximately $1 \sim 2$ million cells to accurately capture the 3D detached shock systems and tip leakage vortices. The evaluation metric used the relative root mean square error (RRMSE) defined by PLAID.
For a physical field $U$, RRMSE is defined as:
\begin{equation}
\mathrm{RRMSE} = \left( \frac{1}{n_*} \sum_{i=1}^{n_*} \frac{\|U_{\mathrm{ref}}^i - U_{\mathrm{pred}}^i\|_2^2}{N \cdot \|U_{\mathrm{ref}}^i\|_\infty^2} \right)^{1/2} .
\end{equation}
For a scalar output $\omega$, RRMSE is defined as:
\begin{equation}
\mathrm{RRMSE}_s = \left( \frac{1}{n_*} \sum_{i=1}^{n_*} \frac{|\omega_{\mathrm{ref}}^i - \omega_{\mathrm{pred}}^i|^2}{|\omega_{\mathrm{ref}}^i|^2} \right)^{1/2} .
\end{equation}

We globally adopted a two-stage optimization strategy, combining "Feature Auto-Encoder (FAE) warm-up" with "Flow Matching dynamics optimization".

\textbf{Baseline Selection:}
We constructed a baseline matrix covering mainstream architectures to evaluate performance boundaries. This selection targeted the strong nonlinearity and multi-scale topology of fluid data. While traditional physics-informed neural networks (PINNs) have achieved tremendous success in many PDEs, they can occasionally encounter soft constraint gradient conflicts and convergence challenges under highly compressible 3D high Mach number flows. Thus, the baseline library focused on stable proxy models. These were divided into three paradigms:
\begin{itemize}
    \item \textbf{Traditional ROMs:} We introduced POD-MLP as a linear benchmark. This verified the advantages of deep implicit manifolds over traditional truncation for strong nonlinear features like transonic shocks.
    \item \textbf{Discriminative Operators:} We included PointNet~\cite{qi2017pointnet}, PointNet++~\cite{qi2017pointnetplusplus}, DGCNN~\cite{wang2019dynamic}, GINO~\cite{li2023geometry}, and Point Transformer for the surface domain. For the volumetric field, we selected DeepONet~\cite{lu2021DeepONet}, 3D-FNO~\cite{li2021fourier}, 3D U-Net~\cite{cicek20163dunet}, and MeshGraphNet (MGN).
    \item \textbf{Generative Architectures:} We also introduced the conditional denoising diffusion probabilistic model (c-DDPM) and 3D conditional variational autoencoder (3D C-VAE)~\cite{sohn2015learning}. This ensured a fair comparison of generative modeling capabilities.
\end{itemize}

\subsection{Few-Shot Aerodynamic Prediction for BlendedNet}

\begin{table}[!htbp]
\caption{Comprehensive comparison of aerodynamic prediction performance and data efficiency on BlendedNet}
\label{tab_blendednet}
\centering
\small 
\begin{tabular}{lcccc}
\hline
Model & Number of Samples & $C_{p}$ MAE & $C_{fx}$ MAE & $C_{fz}$ MAE \\
\hline
\textit{Generative Baselines} & & & & \\
c-DDPM (Diffusion) & 100 & $2.94 \times 10^{-1}$ & $6.72 \times 10^{-3}$ & $4.59 \times 10^{-3}$ \\
\textbf{GeoFunFlow-3D (Ours)} & \textbf{100} & \textbf{$1.43 \times 10^{-1}$} & \textbf{$4.99 \times 10^{-3}$} & \textbf{$1.82 \times 10^{-3}$} \\
\textbf{GeoFunFlow-3D (Ours)} & \textbf{500} & \textbf{$1.04 \times 10^{-1}$} & \textbf{$4.23 \times 10^{-3}$} & \textbf{$1.62 \times 10^{-3}$} \\
\hline
\textit{Discriminative Baselines} & & & & \\
PointNet & 100 & $2.97 \times 10^{-1}$ & $7.81 \times 10^{-3}$ & $4.21 \times 10^{-3}$ \\
GINO & 100 & $2.28 \times 10^{-1}$ & $8.46 \times 10^{-3}$ & $4.83 \times 10^{-3}$ \\
DGCNN & 100 & $2.14 \times 10^{-1}$ & $5.21 \times 10^{-3}$ & $2.38 \times 10^{-3}$ \\
PointNet2 & 100 & $1.99 \times 10^{-1}$ & $5.27 \times 10^{-3}$ & $3.64 \times 10^{-3}$ \\
Point Transformer & 100 & $1.93 \times 10^{-1}$ & $5.12 \times 10^{-3}$ & $3.40 \times 10^{-3}$ \\
\hline
\textit{Oracle (Upper Bound)} & & & & \\
PointNet+FiLM (Official) & 8830 & $3.72 \times 10^{-2}$ & $1.35 \times 10^{-3}$ & $7.96 \times 10^{-4}$ \\
\textbf{GeoFunFlow-3D (FAE)} & \textbf{8830} & \textbf{$1.01 \times 10^{-1}$} & \textbf{$4.01 \times 10^{-3}$} & \textbf{$1.60 \times 10^{-3}$} \\
\hline
\end{tabular}
\end{table}

\begin{figure*}[!htbp]
\centering
\includegraphics[width=0.95\textwidth]{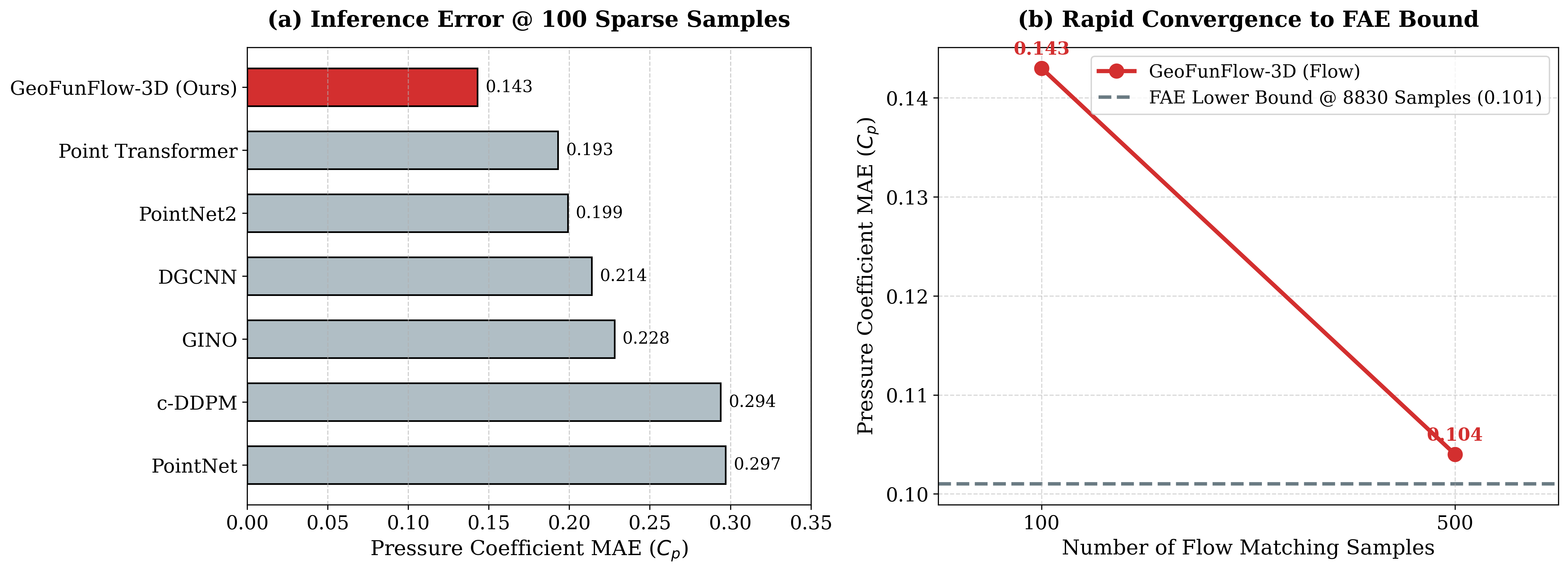} 
\caption{Data efficiency and scaling behavior of GeoFunFlow-3D on the BlendedNet dataset. (a) Comparison of $C_p$ MAE under the extreme low-data regime (100 samples). (b) Rapid error convergence of the flow matching stage towards the theoretical FAE lower bound as sample size increases.}
\label{fig_blendednet_efficiency}
\end{figure*}

We conducted benchmark tests using only 100 and 500 samples to evaluate generalization under restricted data, following the LAMP framework~\cite{nehme2026lamp}. Table~\ref{tab_blendednet} presents the comprehensive performance comparisons across all aerodynamic coefficients ($C_p$, $C_{fx}$, $C_{fz}$) under this few-shot setting. To visually highlight the data efficiency, Figure~\ref{fig_blendednet_efficiency} illustrates the scaling behavior specifically for the pressure coefficient.

Under the extreme 100-sample condition (Figure~\ref{fig_blendednet_efficiency}a), baseline models encountered distinct challenges. For instance, the purely data-driven c-DDPM experienced mode collapse ($C_p$ MAE of $2.94 \times 10^{-1}$), likely because a standard Markov chain lacking explicit physical scheduling often faces inherent challenges in establishing a well-posed score function within a highly sparse probability space. Similarly, while introducing a powerful self-attention mechanism, the Point Transformer's accuracy reached a plateau at $1.93 \times 10^{-1}$ under this sparse setting. In stark contrast, GeoFunFlow-3D demonstrated robust performance, reducing the error to $1.43 \times 10^{-1}$. This stability is theoretically attributed to the minimum transport cost of the approximate $W_2$ geodesic and the spectral boundedness of the No-AD engine.

Furthermore, Figure~\ref{fig_blendednet_efficiency}b reveals a steep data-efficiency scaling law. With an increase to just 500 samples, the inference error of GeoFunFlow-3D drops sharply to $1.04 \times 10^{-1}$, closely approaching the physical lower bound ($1.01 \times 10^{-1}$) established by the Feature Auto-Encoder (FAE) trained on the full dataset (8830 samples). By leveraging the dual optimization dynamics, our framework effectively fitted the complex joint probability distribution while consuming only 5.7\% of the official baseline data volume. This rapid convergence confirms that the proposed physics-guided flow matching formulation effectively overcomes data sparsity.

Objectively, the absolute FAE accuracy does not surpass that of the official PointNet pure coordinate operator. This reflects a "Generative Trade-off". The continuous space must be compressed into a discrete latent grid for unconditional manifold evolution. This quantization smoothing sacrifices some high-frequency limit accuracy. However, it provides global modeling capability for flow field uncertainty.

The proposed framework achieved superior spatial error control compared to the baselines (see visual comparisons in Figure~\ref{app:error_improvement}). This advantage was further confirmed in the generation of the surface friction coefficient ($C_{fx}$).

\begin{figure*}[!htbp]
\centering
\includegraphics[width=0.8\textwidth]{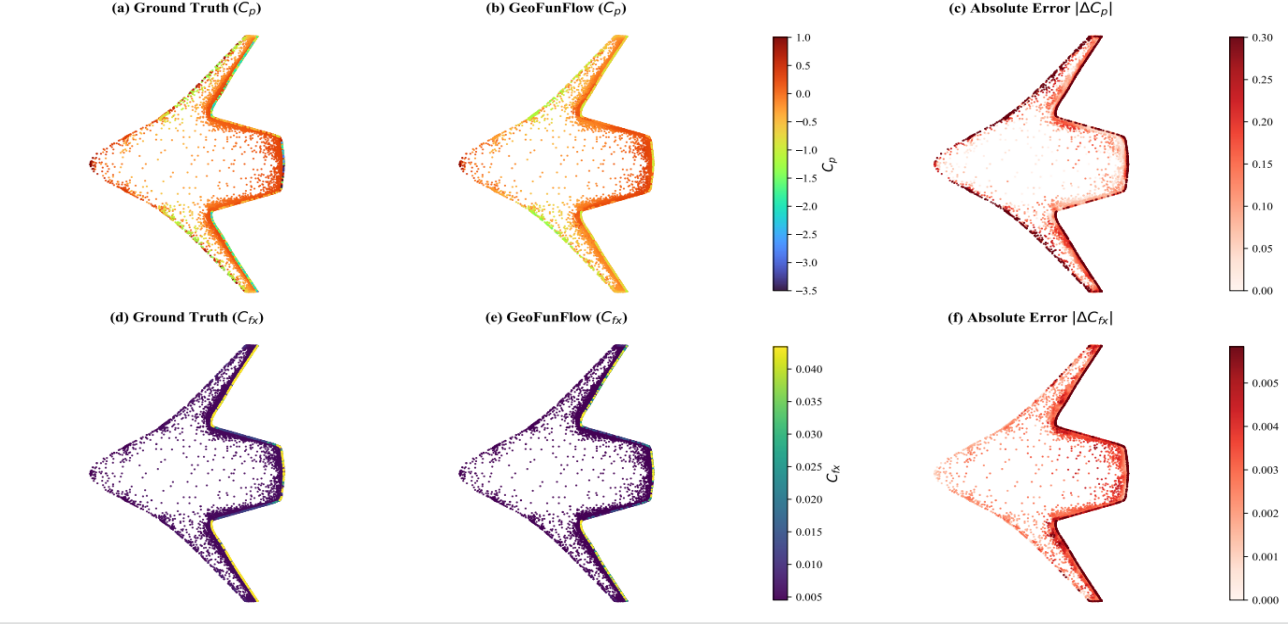} 
\caption{Generation comparison of pressure coefficient ($C_p$) and surface friction coefficient ($C_{fx}$) for the same blended-wing-body.}
\label{fig_cp_cfx}
\end{figure*}

As shown in Figure~\ref{fig_cp_cfx}, the surface friction magnitude is extremely small ($\le 10^{-3}$) and highly sensitive to curvature changes. However, GeoFunFlow-3D successfully captured the micro-friction fluctuations at the wing-body junction and the trailing edge. This was achieved through the high-frequency residual compensation of SATO (neural MsFEM basis functions), demonstrating the framework's capability in decoupling macroscopic and microscopic scales.

Finally, the framework achieved a smooth evolution from Gaussian noise to the physical field (see multi-time-step trajectory in Figure~\ref{app:generation}). The generative probability space also showed strong physical immunity to initial uncertainty, stably converging to similar flow fields even when tested with extreme noise starting points (see convergence tests in Figure~\ref{app:extreme_noise}).

\subsection{Complex Thermodynamic Volume Field Generation for Rotor37}

\begin{table*}[!htbp]
\caption{Comprehensive comparison of architecture representation capability and generalization blind test on Rotor37}
\label{tab_rotor}
\centering
\small 
\begin{tabular}{llcccccc}
\hline
\multirow{2}{*}{Model} & \multirow{2}{*}{Architecture} & \multicolumn{3}{c}{1000 Samples (Train Limit)} & \multicolumn{3}{c}{200 Samples (Blind Test)} \\
\cline{3-5} \cline{6-8}
 & & Pressure & Density & Temperature & Pressure & Density & Temperature \\
\hline
POD-MLP & Traditional ROM & 0.0164 & 0.0165 & 0.0038 & - & - & - \\
3D C-VAE & Traditional Generative & 0.1390 & 0.1450 & 0.0469 & 0.1145 & 0.1208 & 0.0356 \\
DeepONet & Pure Coordinate Operator & 0.0813 & 0.0828 & 0.0140 & 0.1730 & 0.1800 & 0.0426 \\
3D-FNO & Frequency-Domain Operator& 0.0259 & 0.0262 & 0.0119 & 0.0840 & 0.0836 & 0.0086 \\
3D U-Net & Voxel Regression & 0.0253 & 0.0259 & 0.0117 & 0.1206 & 0.1238 & 0.0460 \\
\textbf{GeoFunFlow-3D} & \textbf{Flow Matching Gen.} & \textbf{0.0121} & \textbf{0.0137} & \textbf{0.0065} & \textbf{0.0215} & \textbf{0.0218} & \textbf{0.0134} \\
\hline
MGN & Graph Neural Network & - & - & - & 0.0170 & 0.0114 & 0.0140 \\
GCNN & Graph Neural Network & - & - & - & 0.0170 & 0.0214 & 0.0039 \\
\hline
\end{tabular}
\end{table*}

In the 1000-sample limit test (Table~\ref{tab_rotor}), POD-MLP reconstructed the temperature field well. However, due to its inherently linear nature, it naturally encountered truncation limitations at sharp transonic shock discontinuities. Consequently, its pressure RRMSE was constrained to 0.0164. GeoFunFlow-3D overcame this linear bottleneck via nonlinear Flow Matching evolution. Its pressure error was reduced to 0.0121.

In the 200-sample blind test, while DeepONet and 3D U-Net provide strong baselines, they understandably exhibited higher error margins under such extreme data scarcity. Meanwhile, physical laws (equation of state and entropy limits) constrained GeoFunFlow-3D. Thus, it maintained a low pressure error of 0.0215. Its core metrics matched or exceeded graph neural networks (MGN, GCNN) specialized for unstructured grids.

\begin{figure*}[!htbp]
\centering
\includegraphics[width=0.9\textwidth]{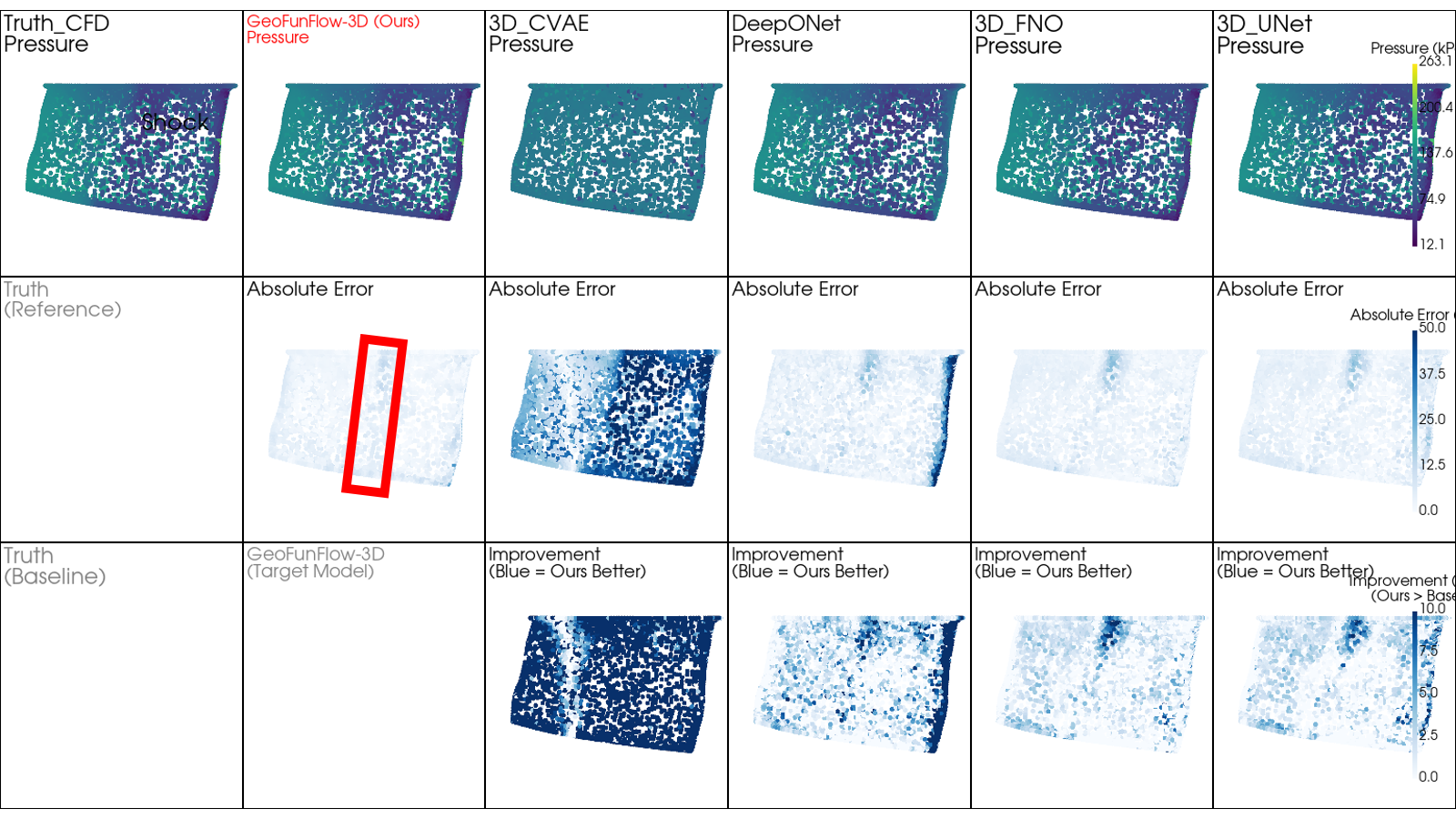} 
\caption{3D visual comparison matrix of the pressure field for the Rotor37 transonic compressor. The middle row shows the absolute error (dark red indicates higher error), while the bottom row highlights accuracy improvements (dark blue indicates higher improvement). Zoom-in panels display localized absolute errors near the shock region.}
\label{fig_rotor37_matrix}
\end{figure*}

The 3D spatial error matrix (Figure~\ref{fig_rotor37_matrix}) reveals reconstruction details for the high Mach number internal flow. Significant transonic detached shocks exist in the middle and rear blade sections of the true flow field (first column, \textit{Truth\_CFD}). A red arrow indicates this topology, while a solid red box in the reference error map outlines the localized zoom-in region. Faced with this nonlinear topology, the 3D-CVAE baseline experienced inherent variational lower bound smoothing, making it challenging to maintain sharp shock fronts. Without explicit physical shock constraints, pioneering conventional operators (such as DeepONet and 3D-FNO) naturally exhibited some localized error accumulation near the shock surface. 

In contrast, a local entropy increase penalty was explicitly injected into GeoFunFlow-3D. As a result, physical discontinuities were effectively restored. A more uniform full-field error distribution was achieved. The third row of the matrix shows dark blue regions representing significant accuracy improvements. These regions spatially overlap with the marked 3D shock topology. This visually confirms the effectiveness of volumetric thermodynamic constraints in suppressing high-frequency mutation errors.

Overall, GeoFunFlow-3D demonstrated competitive spatial fidelity in generating the transonic volume field (\textbf{detailed slice comparisons for density and pressure are shown in Figure~\ref{app:rotor_fields}}).

\subsection{Inference Efficiency Analysis for Engineering Deployment}

\begin{table*}[!htbp] 
\caption{Comprehensive comparison of model complexity and inference acceleration strategies}
\label{tab_efficiency}
\centering
\small 
\begin{tabular}{llccp{4cm}} 
\hline
Dataset & Model & Parameters (M) & Single-Sample Time (ms) & Remarks \\
\hline
\multirow{8}{*}{BlendedNet} 
 & PointNet & 0.88 & 2.6 & \\
 & DGCNN & 0.97 & 284.2 & \\
 & GINO & 0.10 & 24.7 & \\
 & PointNet2 & 33.63 & 32.6 & \\
 & Point Transformer & 0.17 & 19.1 & Configured for few-shot limits \\
 & c-DDPM & 0.13 & 387.2 & 1000-step reverse sampling \\
 & \textbf{GeoFunFlow-3D (Base)} & \textbf{37.36} & \textbf{3590.60} & \\
 & \textbf{GeoFunFlow-3D (Cached)}& \textbf{37.36} & \textbf{2645.48} & \textbf{26.3\% Speedup} \\
\hline
\multirow{7}{*}{Rotor37} 
 & POD-MLP & 0.04 & 0.1 & Linear ROM \\
 & DeepONet & 0.13 & 0.69 & \\
 & 3D U-Net & 5.61 & 2.17 & \\
 & 3D C-VAE & 6.21 & 4.30 & \\
 & 3D-FNO & 33.58 & 3.31 & \\
 & \textbf{GeoFunFlow-3D (Base)} & \textbf{494.61} & \textbf{17913.06} & \\
 & \textbf{GeoFunFlow-3D (Shared)}& \textbf{154.82} & \textbf{4909.31} & \textbf{72.6\% Speedup} \\
\hline
\end{tabular}
\end{table*}

Engineering deployment requires high inference efficiency. As shown in Table~\ref{tab_efficiency}, Flow Matching models using ODE evolution are naturally slower than direct mapping operators. To resolve this, specific acceleration strategies were designed. We introduced Cross-step Feature Caching in BlendedNet, achieving a 26.3\% speedup, and applied Cross-layer Weight Sharing in Rotor37. This reduced redundant parameters by 68.7\% and accelerated single-sample generation by 72.6\%. The flow matching model needs multi-step ODE integration. Thus, the single-sample inference takes about 4.9 seconds. This is slower than single-step networks. However, this speed represents a substantial improvement for computational physics. For the same mesh with one to two million cells, a traditional commercial RANS solver (e.g., Ansys Fluent or Siemens STAR-CCM+) typically requires from tens of minutes to several hours depending on the available high-performance computing resources to reach equivalent convergence. Therefore, the proposed framework achieves a speedup of over $10^3$ times. Compared to computationally intensive traditional CFD solvers, this optimized inference time aligns well with engineering acceleration requirements. 

Furthermore, it is important to objectively acknowledge the inherent generative trade-off associated with the proposed framework: the ultra-fast online inference speed is amortized by a substantial offline training cost. For transparency, the offline training of the full GeoFunFlow-3D model (e.g., 494.61M parameters for the Rotor37 volumetric field) required approximately 72 to 144 hours on a computation node equipped with 2 NVIDIA Tesla V100 (32GB) GPUs to reach convergence. However, in practical engineering design loops where thousands of aerodynamic evaluations are required, this one-time intensive offline cost is highly acceptable. It successfully shifts the heavy computational burden of traditional CFD solvers into the offline phase, enabling near real-time, high-fidelity 3D volume field generation during the online inference phase.

\subsection{Core Module Ablation Study}

\begin{table}[!htbp]
\caption{Ablation Study of core modules in GeoFunFlow-3D}
\label{tab_ablation}
\centering
\small
\begin{tabular}{lccp{5cm}}
\hline
Model & $C_{p}$ MAE & $C_{fx}$ MAE & Remarks \\
\hline
\textbf{GeoFunFlow-3D} & \textbf{$1.43 \times 10^{-1}$} & \textbf{$4.99 \times 10^{-3}$} & \textbf{Full model} \\
w/o SATO & $1.76 \times 10^{0}$ & $4.70 \times 10^{-2}$ & Removed topological refinement decoder \\
w/o TV Loss & $2.79 \times 10^{-1}$ & $6.04 \times 10^{-3}$ & Removed total variation regularization \\
w/o Hard Mask & $2.95 \times 10^{-1}$ & $6.23 \times 10^{-3}$ & Removed deterministic geometric mask \\
PointNet-Flow & $2.97 \times 10^{-1}$ & $7.81 \times 10^{-3}$ & Removed latent space and DiT \\
\hline
\end{tabular}
\end{table}

\begin{figure*}[!htbp]
\centering
\includegraphics[width=0.6\textwidth]{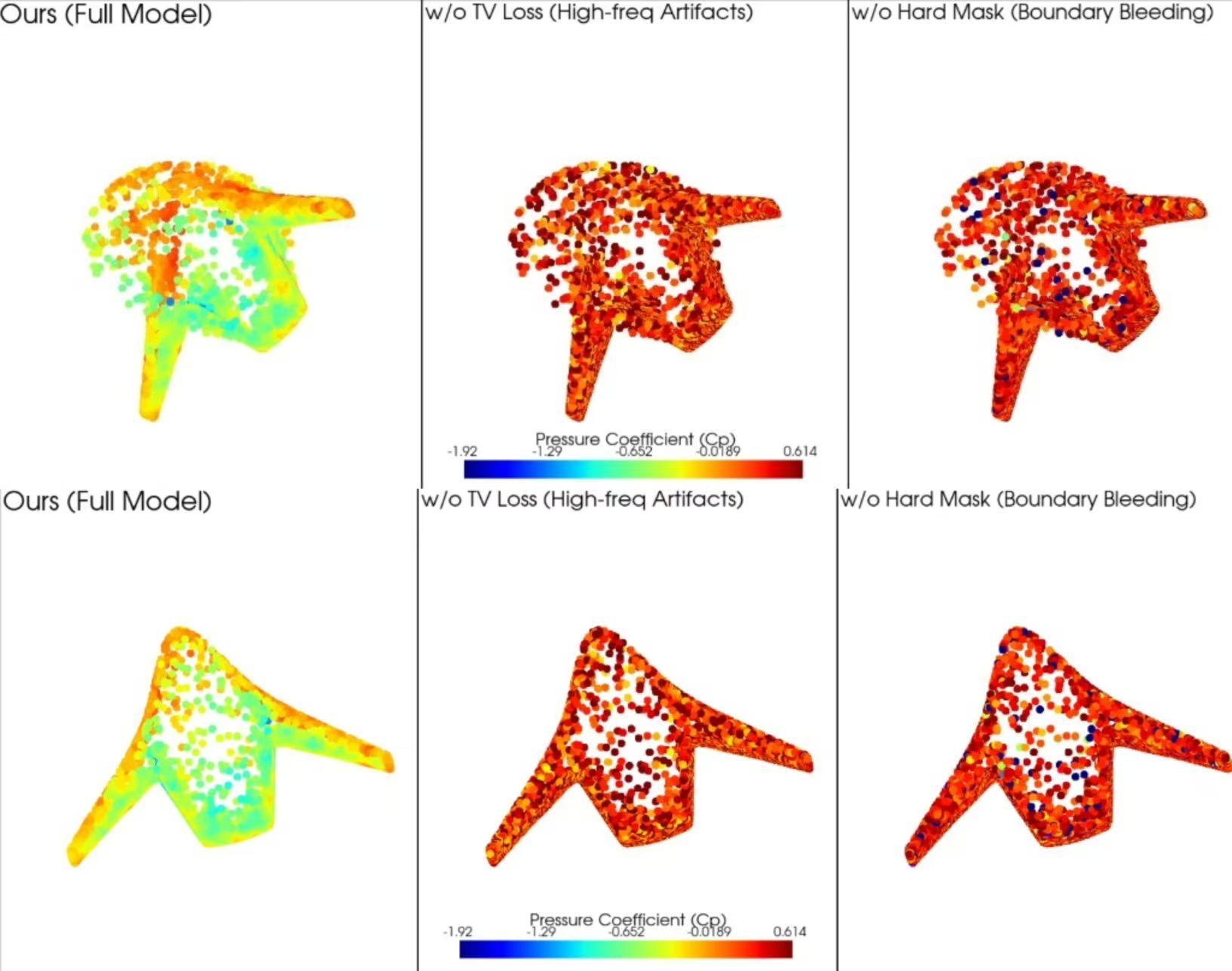} 
\caption{Impact of missing core physical operators on flow field generation. Dark patches indicate numerical anomaly regions deviating from physical norms.}
\label{fig_ablation_visual}
\end{figure*}

The Ablation Study (Table~\ref{tab_ablation} and Figure~\ref{fig_ablation_visual}) demonstrates the specific contributions of the mathematical mechanisms:
\begin{itemize}
    \item \textbf{High-frequency oscillation without TV Loss:} We removed the total variation functional of the inverse Jacobian pullback. Consequently, numerous abnormal extreme points appear in the generated field. This represents the Gibbs phenomenon caused by lacking a flux limiter. It highlights the necessity of anisotropic smoothing for spatial dynamic stability.
    \item \textbf{Cross-boundary leakage without Hard Mask:} We removed the implicit phase field mask. Grid nodes in the non-fluid domain unrestrictedly participated in gradient optimization. This allowed erroneous energy to spill out, triggering cross-boundary leakage in the external flow field. This indicates that geometric hard constraints are highly beneficial in complex open domains.
\end{itemize}

\subsection{Uncertainty Quantification and Physics Compliance}

\begin{figure}[!htbp]
\centering
\includegraphics[width=1.0\textwidth]{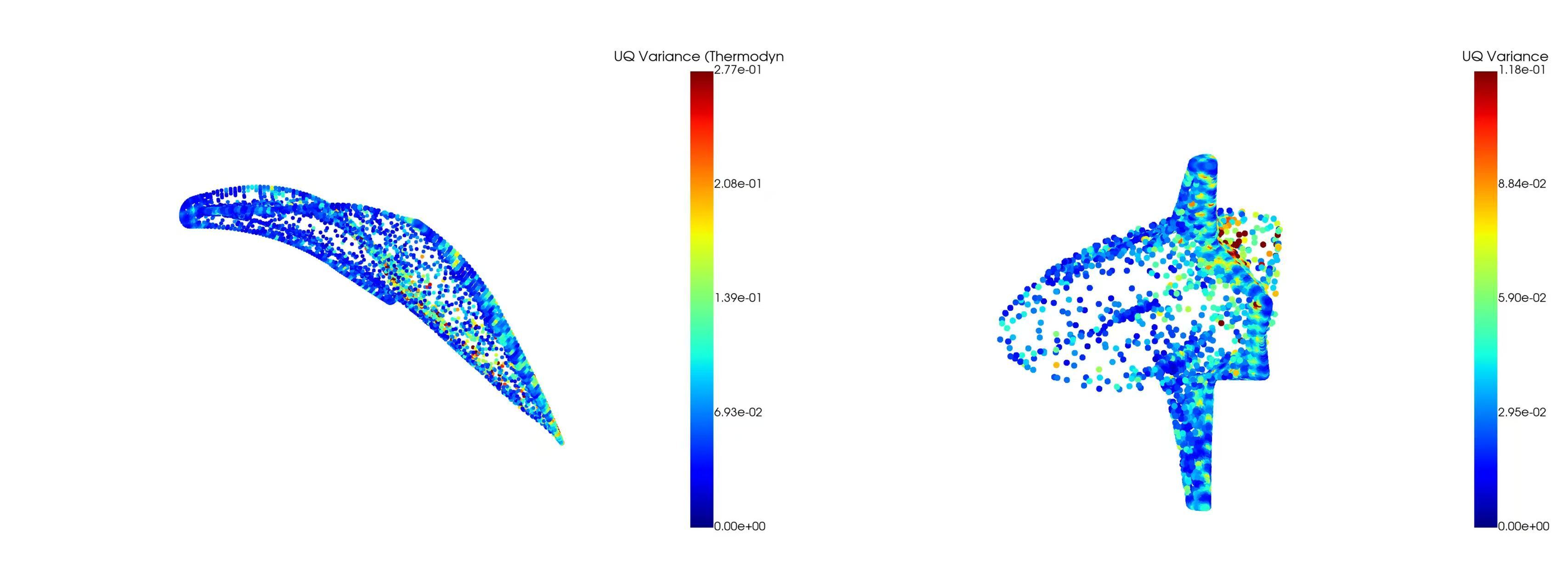} 
\caption{Spatial Uncertainty Quantification (UQ) variance distribution. The upper figure shows the Rotor37 transonic internal flow. The lower figure shows the BlendedNet external aerodynamic surface.}
\label{fig_combined_uq}
\end{figure}

To evaluate result confidence, we injected 30 independent Gaussian noise sets into the latent space and calculated the spatial variance distribution (UQ) of the generated fields. As shown in Figure~\ref{fig_combined_uq}, the model exhibited polymorphic variance adaptation. For the transonic internal flow (upper figure), the high-variance region accurately matched the detached shock front. For the external aerodynamic field (lower figure), uncertainty concentrated in the wing-body interference and wake separation regions. In laminar smooth regions, the generated trajectory showed high certainty with near-zero variance. This indicates that the model accurately maps the intrinsic physical chaos of the Navier-Stokes system using probability distributions.

\begin{figure}[!htbp]
\centering
\includegraphics[width=0.6\textwidth]{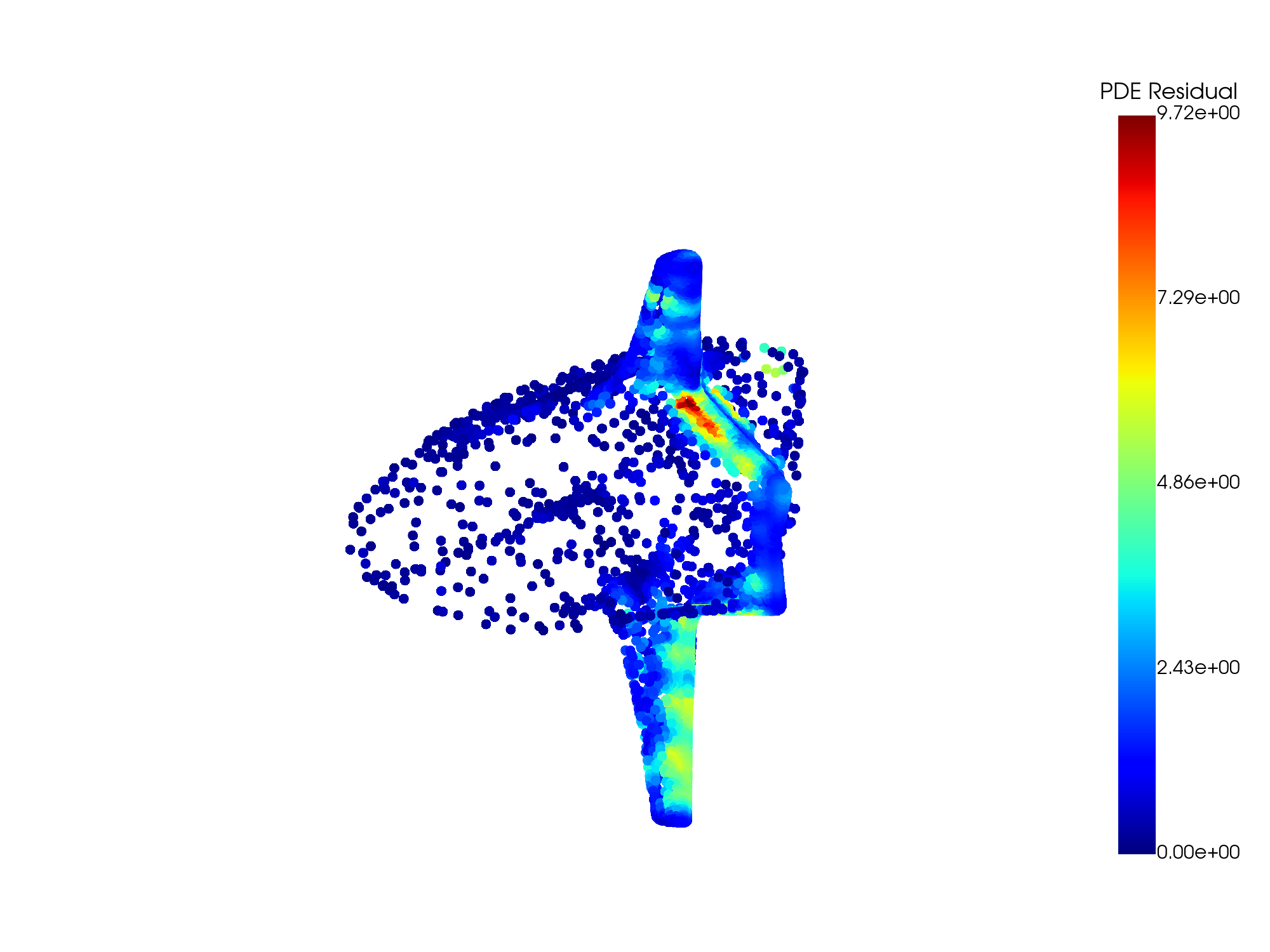} 
\caption{Spatial distribution of PDE physical residuals of the generated flow field calculated based on the discrete differential engine.}
\label{fig_blendednet_pde}
\end{figure}

Finally, we substituted the predicted mean flow field back into Euclidean space and calculated the physical partial differential equation (PDE) residuals (Figure~\ref{fig_blendednet_pde}). The generated flow field maintained extremely low physical residuals (dark blue) across most smooth computational domains. Local residual peaks (red highlights) did not diverge. Instead, they consistently localized at the wing-body junction and trailing edge regions. This aligns well with fluid mechanics knowledge. Strong geometric curvature mutations and viscous flow separation exist here. These cause continuous differential operators to generate truncation error peaks. Overall, this residual distribution suggests that GeoFunFlow-3D maintains a high degree of consistency with global mass and momentum conservation laws.

\begin{figure}[!htbp]
\centering
\includegraphics[width=0.6\textwidth]{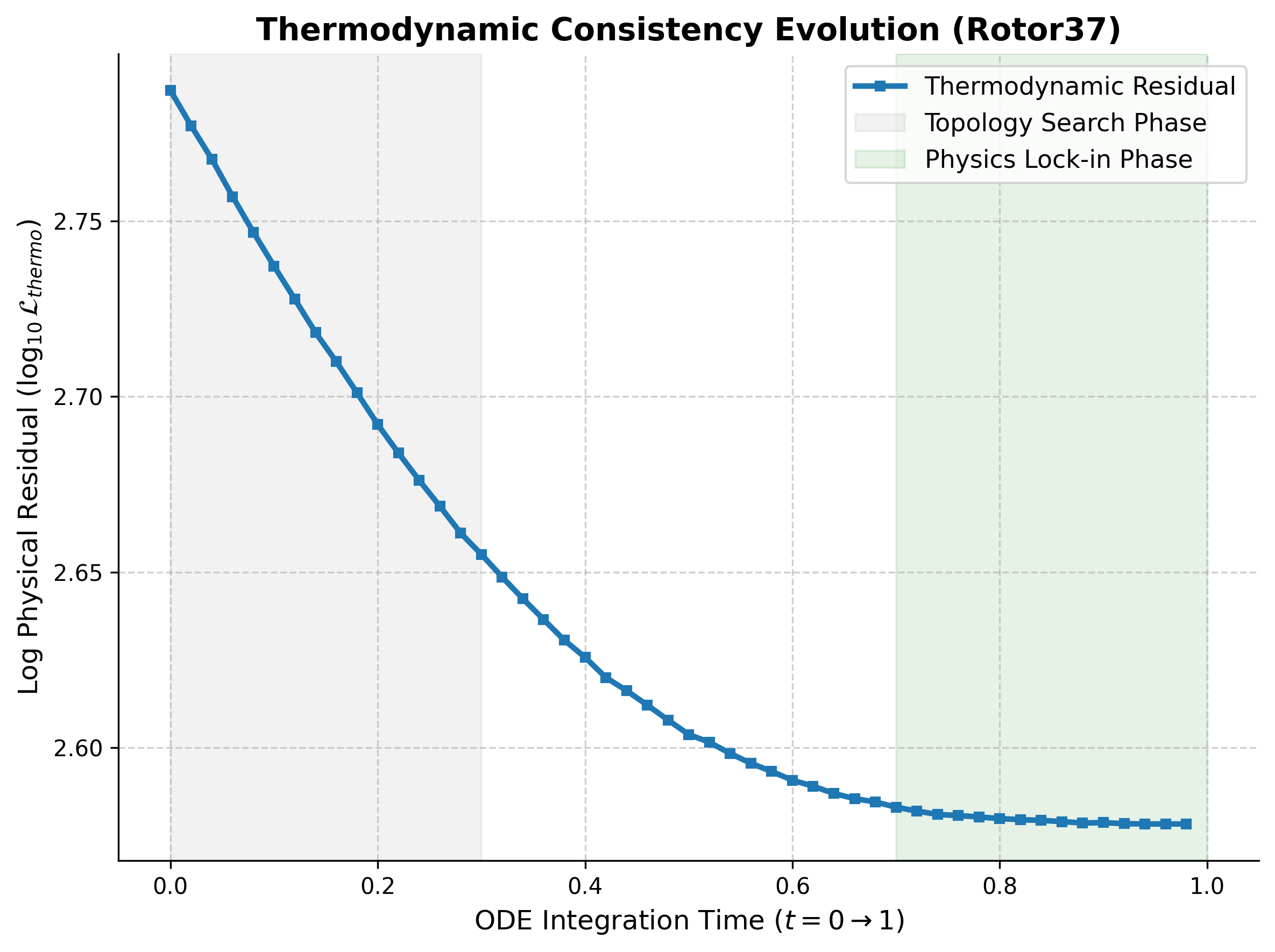} 
\caption{Time evolution curve of the global thermodynamic residual during the Flow Matching ODE generation trajectory ($t=0 \to 1$). The figure divides the pure data-driven "Topology Search Phase" and the constraint-dominated "Physics Lock-in Phase".}
\label{fig_residual_convergence}
\end{figure}

To further explore dynamic physical consistency during generation, we tracked the global volumetric thermodynamic residual distribution and monitored the model evolution along the ODE generation geodesic using the Rotor37 dataset (Figure~\ref{fig_residual_convergence}). The horizontal axis represents the flow matching ODE integration time ($t=0 \to 1$). The vertical axis represents the mean physical residual ($\mathcal{L}_{thermo}$) plotted on a logarithmic scale. The legend indicators include the "Thermodynamic Residual" curve, the "Topology Search Phase" (gray background, $t \in [0, 0.3]$), and the "Physics Lock-in Phase" (green background, $t > 0.7$). The latter two are driven by the dual optimization dynamics scheduling.

The Rotor37 dataset contains complex transonic features with strict physical constraints. During the initial Topology Search Phase, the system started from pure Gaussian white noise. The flow field was highly chaotic. The model prioritized the spatial alignment of macroscopic density. Variational Homotopy Scheduling dynamically suppressed weights. This protected the model from divergence risks caused by high initial residuals. Consequently, the residual curve exhibited a steep monotonic decline. During the Physics Lock-in Phase, the macroscopic aerodynamic topology was preliminarily established. Physical constraints began to dominate the precise calibration of microscopic discontinuous structures like detached shocks. The residual curve transitioned into a stable exponential decay and asymptotically locked at a minimum level. This temporal convergence process indicated that the framework effectively coordinates data generation and physical conservation laws. It ensures the generated 3D volumetric flow field maintains strong physical consistency.

\section{Conclusions}
In this study, we proposed GeoFunFlow-3D, a generative surrogate model for 3D aerodynamic and thermodynamic fields. This framework deeply integrates the No-AD high-order discrete operator, the Lagrangian phase-field hard mask, and the flow matching network using Variational Homotopy Scheduling. It systematically mitigates multi-scale errors and high-frequency oscillations in 3D unstructured flow field prediction. This framework achieves reliable inference of complex physical fields under highly restricted sample sizes.

\textbf{Limitations and Future Outlook:} GeoFunFlow-3D shows effectiveness on conventional continuous curvature mutations. However, for 3D geometries with multi-scale topological defects (e.g., multiple holes, extreme aspect ratios), the inverse Jacobian matrix of the global homeomorphic mapping may still experience local singularity. This remains a limitation for coordinate-mapping neural operators. Additionally, the pure geometry-prior generation depends on the latent manifold. The FAE warm-up phase builds this manifold using source domain data. The model achieves zero-shot inference within this training distribution. However, this generalization capability has strict boundaries. For entirely new and out-of-distribution geometries, the model may experience performance degradation. In these cases, a small amount of data is needed to fine-tune the manifold. This step is intended to support accuracy recovery for highly novel shapes. Future research will deepen the underlying mathematical framework from two dimensions. Algebraically, we will explore neural operator preconditioning based on the Jacobian-Free Newton-Krylov (JFNK) subspace method to reshape the characteristic spectrum of ill-conditioned matrices. Functionally, we will construct weak-form reconstruction based on local energy functionals. This draws inspiration from atlas-based overset mapping theory and the Deep Ritz method~\cite{weinan2018deep}. Furthermore, while the current framework achieves high-fidelity point-wise reconstruction of microscopic field variables, it has not yet explicitly integrated macroscopic aerodynamic and thermodynamic integral metrics (e.g., total pressure recovery ratio and adiabatic efficiency for compressor rotors) into the evaluation loop. Incorporating these macroscopic engineering quantities as volume-integral soft constraints into the generative loss function will be a crucial next step to fully align with turbomachinery design requirements. The spectral-spatial-temporal joint physical constraint paradigm established by GeoFunFlow-3D provides a solid theoretical foundation for the intelligent forward design of next-generation aerospace vehicles, demonstrating immense potential for enabling rapid aerodynamic shape optimization loops in complex engineering scenarios.

\section*{Data Availability}
\indent Data will be made available on request.

\section*{Acknowledgements} 
\emph{The authors are very much indebted to the referees for their constructive comments and valuable suggestions}.

\clearpage
\appendix

\section{Supplementary 3D Visualizations of the Generation Pipeline} \label{App:3D_Pipeline}

\begin{figure}[htbp]
    \centering
    \includegraphics[width=0.9\linewidth]{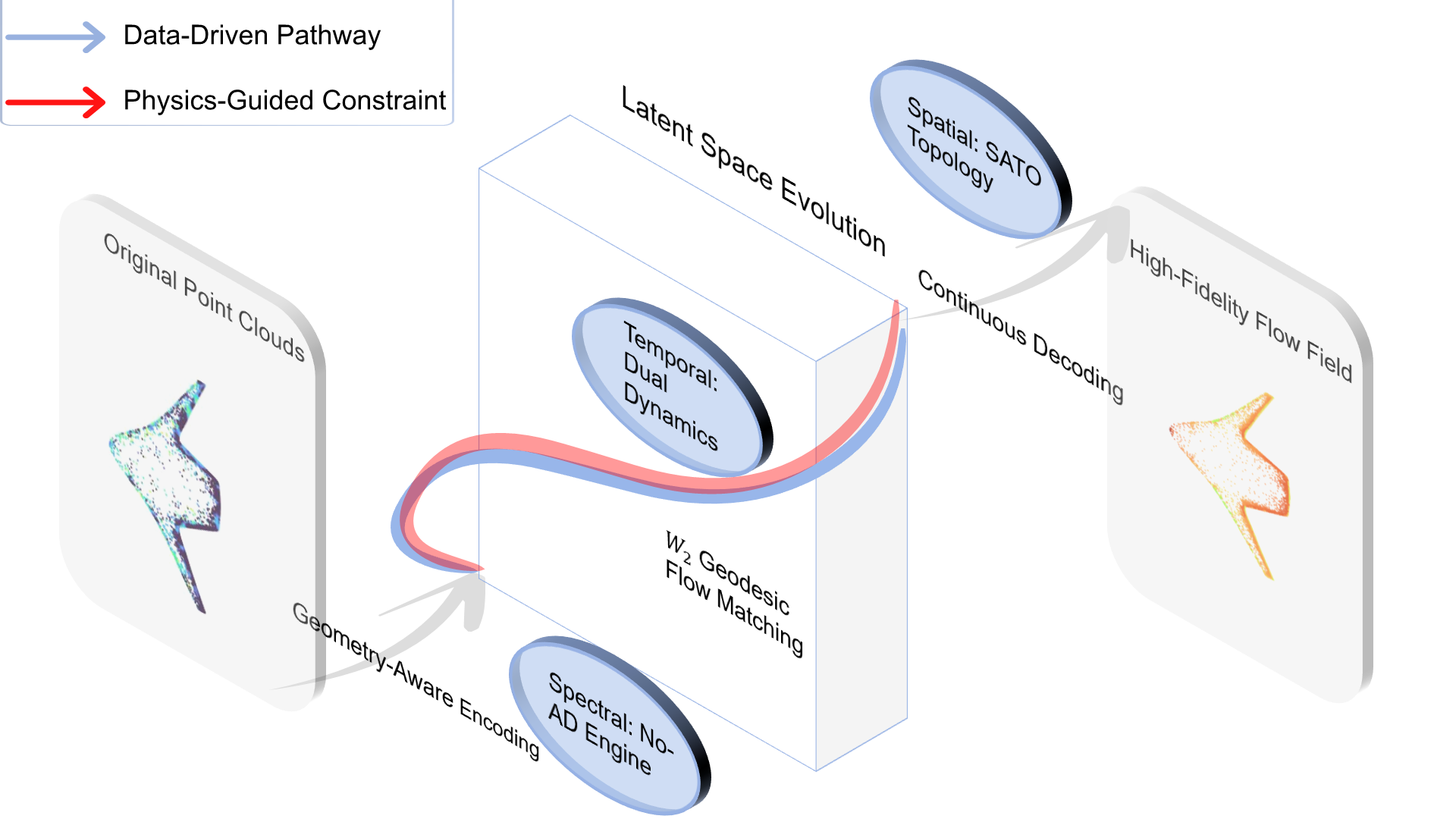}
    \caption{A 3D perspective visualizing the complete zero-observation forward generation pipeline of GeoFunFlow-3D, transitioning from original point clouds to high-fidelity aerodynamic fields.}
    \label{fig:3d_architecture}
\end{figure}

To provide an intuitive spatial perspective of the zero-observation forward generation paradigm discussed in Section 4.1, Figure~\ref{fig:3d_architecture} visualizes the complete pipeline in a 3D context. This explicitly demonstrates the transition from pure geometric point clouds to high-fidelity flow fields within the latent space.\\

\clearpage
\section{Proof of Lipschitz Continuity and Density Invariance of the NW Projection Operator} \label{App:Lipschitz}

\subsection{Proof of Lipschitz Continuity}
In Section 3.1.3, the data format conversion is completed at the engineering level. Furthermore, the underlying algebraic structure provides a theoretical guarantee. 
\begin{proposition}[Lipschitz Continuity]
The physical features can be uniformly and continuously mapped onto the latent tensor grid $\tilde{\phi}$. This holds regardless of significant distortions in the aircraft's external topology or high variance in the point cloud density.
\end{proposition}

\begin{proof}
The projection function is defined as $f(\xi) = \tilde{\phi}(\xi)$. To prove its smoothness under coordinate variations in the latent space, its derivative with respect to the coordinate $\xi$ is examined.
Based on the structure of the Nadaraya-Watson estimator, the normalized kernel weight field is defined as:
\begin{equation}
W_j(\xi) = \frac{\omega(\|\xi - x_j\|)}{\sum_{m} \omega(\|\xi - x_m\|)} .
\end{equation}
The inverse distance kernel $\omega(d) = (d + \epsilon)^{-1}$ selected in this framework satisfies $C^\infty$ continuity when $\epsilon > 0$. Moreover, its denominator (the sum of weights) is strictly greater than zero when $\mathcal{N} \neq \emptyset$. Therefore, $W_j(\xi)$ is non-singular within the domain.

For a tiny displacement $\Delta \xi$, both $\mathcal{K}_\phi$ and $W_j$ are composed of smooth functions. Consequently, their first-order partial derivatives are bounded:
\begin{equation}
\begin{aligned}
\|\nabla_\xi \tilde{\phi}\| &\le \sum_j \left( \|\nabla_\xi W_j\| \cdot \| \text{Feature}_j \| + \|W_j\| \cdot \|\nabla_\xi \mathcal{K}_\phi\| \right) \\
&\le L .
\end{aligned}
\end{equation}
From this, it can be obtained:
\begin{equation}
\|\tilde{\phi}(\xi_1) - \tilde{\phi}(\xi_2)\| \le L \|\xi_1 - \xi_2\| .
\end{equation}
\end{proof}
This property ensures a strictly smooth response of the latent space grid under spatial perturbations. It effectively suppresses the generation of spatial discontinuities.

\subsection{Algebraic Explanation of Density Invariance}
Consider two sets of point clouds, $\mathcal{P}_{dense}$ and $\mathcal{P}_{sparse}$, that represent the same geometry. In a local region, the point density of $\mathcal{P}_{dense}$ is assumed to be $M$ times that of $\mathcal{P}_{sparse}$.
If normalization is not performed, then $\tilde{\phi}_{dense} \approx M \cdot \tilde{\phi}_{sparse}$. This causes the feature energy to explode with resolution. By introducing the partition of unity term $\sum \omega(d)$:
\begin{equation}
\tilde{\phi} = \frac{\sum w_j \cdot f_j}{\sum w_j} .
\end{equation}
As the sampling density increases, both the numerator and the denominator grow linearly and synchronously. Their ratio (the convex combination coefficient) remains stable. Topologically, this guarantees that the feature projection depends only on the geometric envelope, rather than the sampling discreteness.

According to the continuity theorem of the Nadaraya-Watson estimator, the selected distance kernel function $\omega(d) = \frac{1}{d + \epsilon}$ is non-singular and continuous everywhere within the domain. Therefore, the weight field generated by its normalization is theoretically continuous in space. The physical significance of this is as follows. When the latent space grid coordinate $\xi_k$ undergoes a tiny movement in space, its assigned inverse distance weight and the extracted translational features will only experience a smooth and gradual change. Abrupt discontinuities are mathematically mitigated. In mathematical analysis, this property, where a tiny change in input only causes a bounded and smooth change in output, is strictly defined as Lipschitz Continuous.

In summary, the convex combination (partition of unity) of the NW kernel regression guarantees the immunity of feature amplitudes to grid distortion. Meanwhile, the continuity of the kernel function guarantees the strict smoothness of the spatial reconstruction gradient. The superposition of these two major mathematical properties proves a key point. The latent space grid $\tilde{\phi}$ constructed based on this formula achieves effective topological dimension reduction for real 3D unstructured flow fields. This significantly reduces the computational overhead of the downstream generative model (DiT) when processing complex geometries.

\clearpage
\section{Integral Derivation of Generation Error Evolution Based on Grönwall's Inequality} \label{App:Gronwall}

\begin{theorem}[Global Generation Error Bound based on Flow Matching]
Under the optimal transport effective coupling, the $2$-Wasserstein distance between the generated distribution $\mathbb{P}_{gen}$ and the true data distribution $\mathbb{P}_{data}$ in the latent space flow matching is upper-bounded by:
\begin{equation}
W_2(\mathbb{P}_{gen}, \mathbb{P}_{data}) \le \mathcal{O}\left( N^{-\frac{1}{d}} \right) \cdot e^L .
\end{equation}
\end{theorem}

\begin{proof}
The generation process of Flow Matching is essentially completed by solving ordinary differential equations (ODEs). The ideal trajectory is denoted as $Z_t$, and the trajectory actually generated by the neural network is denoted as $\hat{Z}_t$. They satisfy the following differential dynamic equations, respectively:
\begin{equation}
\frac{dZ_t}{dt} = v^*(Z_t, t), \quad \frac{d\hat{Z}_t}{dt} = \hat{v}_\theta(\hat{Z}_t, t) .
\end{equation}

For any time $t$ during the evolution process, the upper limit of the error evolution rate between the two trajectories is determined by the direct difference of the instantaneous velocity fields. This is based on the inequality properties of absolute norm differentiation:
\begin{equation}
\frac{d}{dt} \|Z_t - \hat{Z}_t\| \le \left\| \frac{dZ_t}{dt} - \frac{d\hat{Z}_t}{dt} \right\| = \|v^*(Z_t, t) - \hat{v}_\theta(\hat{Z}_t, t)\| .
\end{equation}

An intermediate term $v^*(\hat{Z}_t, t)$ is introduced on the right side. This aims to separate the pure fitting error of the network and the inherent position deviation error. The above equation is analytically bounded using the triangle inequality and the $L$-Lipschitz continuity of the ideal velocity field $v^*$ (i.e., $\|v^*(x) - v^*(y)\| \le L\|x - y\|$):
\begin{equation}
\begin{aligned}
\frac{d}{dt} \|Z_t - \hat{Z}_t\| &\le L \|Z_t - \hat{Z}_t\| + \left\| v^*(\hat{Z}_t, t) - \hat{v}_\theta(\hat{Z}_t, t) \right\| .
\end{aligned}
\end{equation}

The instantaneous trajectory error is defined as $e(t) = \|Z_t - \hat{Z}_t\|$. The upper bound of the functional fitting error is defined as $\epsilon_v(t)$. The above equation is simplified into a linear differential inequality:
\begin{equation}
e'(t) \le L e(t) + \epsilon_v(t) .
\end{equation}

To solve this inequality, both sides are multiplied by an integrating factor $e^{-Lt}$. This constructs a complete derivative form:
\begin{equation}
\begin{aligned}
e^{-Lt}e'(t) - L e^{-Lt} e(t) &\le \epsilon_v(t) e^{-Lt} \\
\implies \frac{d}{dt} \left( e^{-Lt} e(t) \right) &\le \epsilon_v(t) e^{-Lt} .
\end{aligned}
\end{equation}

A definite integral is performed on the above equation from time $t=0$ to $t=1$. All generation trajectories start from the same prior noise point. Therefore, the initial state error is zero (i.e., $e(0) = 0$). The integral result is:
\begin{equation}
e^{-L} e(1) \le \int_0^1 \epsilon_v(t) e^{-Lt} dt .
\end{equation}

Both sides of the equation are multiplied by $e^L$. It is noted that $e^{L(1-t)} \le e^L$ always holds during the evolution time $t \in [0,1]$. Thus, the strict individual error upper bound of the final generated state can be derived:
\begin{equation}
\begin{aligned}
\|Z_1 - \hat{Z}_1\| &\le \int_0^1 \epsilon_v(t) e^{L(1-t)} dt \\
&\le \left( \int_0^1 \epsilon_v(t) dt \right) e^L = \epsilon_v \cdot e^L .
\end{aligned}
\end{equation}

Furthermore, this is generalized to the macroscopic probability measure space. In optimal transport theory, the $W_2(\mathbb{P}_{gen}, \mathbb{P}_{data})$ distance is bounded by the expectation of the squared trajectory error under effective coupling:
\begin{equation}
W_2(\mathbb{P}_{gen}, \mathbb{P}_{data}) \le \sqrt{ \mathbb{E} \left[ \|Z_1 - \hat{Z}_1\|^2 \right] } .
\end{equation}

The previously obtained expected fitting error limit $\epsilon_v \le \mathcal{O}(N^{-\frac{1}{d}})$ due to manifold dimension reduction is substituted into the above equation. Finally, the theoretical upper bound of the global generation error for latent space flow matching is obtained:
\begin{equation}
W_2(\mathbb{P}_{gen}, \mathbb{P}_{data}) \le \mathcal{O}\left( N^{-\frac{1}{d}} \right) \cdot e^L .
\end{equation}
\end{proof}

\clearpage
\section{Derivation of the Fourth-Order Central Difference Discrete Operator Based on Taylor's Expansion} \label{App:Taylor}

GeoFunFlow-3D must capture physical discontinuities such as shocks with high fidelity. It also needs to suppress numerical dissipation caused by truncation errors. To achieve this, a difference operator with a truncation error of $\mathcal{O}(h^4)$ is rigorously derived. This derivation utilizes the method of undetermined coefficients and Taylor's Theorem.

Under a uniform latent grid step $h$, the field variable $\tilde{\phi}$ at point $\xi$ is expanded by two steps to the left and right ($\pm h$ and $\pm 2h$):
\begin{equation}
\begin{aligned}
\tilde{\phi}(\xi + h) &= \tilde{\phi}(\xi) + h\tilde{\phi}' + \frac{h^2}{2}\tilde{\phi}'' + \frac{h^3}{6}\tilde{\phi}''' \\
&\quad + \frac{h^4}{24}\tilde{\phi}^{(4)} + \frac{h^5}{120}\tilde{\phi}^{(5)} + \mathcal{O}(h^6) .
\end{aligned}
\label{eq:taylor_p1}
\end{equation}
\begin{equation}
\begin{aligned}
\tilde{\phi}(\xi - h) &= \tilde{\phi}(\xi) - h\tilde{\phi}' + \frac{h^2}{2}\tilde{\phi}'' - \frac{h^3}{6}\tilde{\phi}''' \\
&\quad + \frac{h^4}{24}\tilde{\phi}^{(4)} - \frac{h^5}{120}\tilde{\phi}^{(5)} + \mathcal{O}(h^6)
\end{aligned}
\label{eq:taylor_m1}
\end{equation}

Equation \eqref{eq:taylor_m1} is subtracted from Equation \eqref{eq:taylor_p1} to eliminate the even-order terms:
\begin{equation}
\begin{aligned}
\tilde{\phi}(\xi + h) - \tilde{\phi}(\xi - h) &= 2h\tilde{\phi}' + \frac{h^3}{3}\tilde{\phi}''' \\
&\quad + \frac{h^5}{60}\tilde{\phi}^{(5)} + \mathcal{O}(h^7) \quad (\text{Equation A}) .
\end{aligned}
\end{equation}

If a conventional second-order difference is divided by $2h$, the main error term remains $\frac{h^2}{6}\tilde{\phi}''' = \mathcal{O}(h^2)$. To eliminate the third-order derivative term, the $\pm 2h$ expansions are introduced and subtracted:
\begin{equation}
\begin{aligned}
\tilde{\phi}(\xi + 2h) - \tilde{\phi}(\xi - 2h) &= 4h\tilde{\phi}' + \frac{8h^3}{3}\tilde{\phi}''' \\
&\quad + \frac{32h^5}{60}\tilde{\phi}^{(5)} + \mathcal{O}(h^7) \quad (\text{Equation B}) .
\end{aligned}
\end{equation}

Through algebraic reorganization, $8 \times (\text{Equation A}) - (\text{Equation B})$ is calculated. This algebraically eliminates the nonlinear term with the coefficient $\tilde{\phi}'''$:
\begin{equation}
\begin{aligned}
&8[\tilde{\phi}(\xi + h) - \tilde{\phi}(\xi - h)] - [\tilde{\phi}(\xi + 2h) - \tilde{\phi}(\xi - 2h)] \\
&\quad = 12h\tilde{\phi}' - \frac{24h^5}{60}\tilde{\phi}^{(5)} + \mathcal{O}(h^7) .
\end{aligned}
\end{equation}

Both sides are divided by $12h$. The fifth-order infinitesimal component is ignored. Thus, the high-fidelity central difference operator with a truncation error of $\mathcal{O}(h^4)$ is obtained:
\begin{equation}
\begin{aligned}
\tilde{\phi}'(\xi) &= \frac{1}{12h} \Big[ -\tilde{\phi}(\xi + 2h) + 8\tilde{\phi}(\xi + h) \\
&\quad - 8\tilde{\phi}(\xi - h) + \tilde{\phi}(\xi - 2h) \Big] + \mathcal{O}(h^4) .
\end{aligned}
\end{equation}
This purely analytical operator is solidified into a static 3D tensor convolution kernel within the network. It enables efficient estimation of multi-order partial derivatives while avoiding the substantial VRAM overhead associated with continuous automatic differentiation.

\clearpage
\section{Supplementary Experimental Results}

This appendix provides supplementary visual analyses of the GeoFunFlow-3D framework on the primary test datasets to further support the theoretical discussions in the main text. We cover model error control, dynamic convergence robustness, and multi-field coupled generation capabilities.

\subsection{Supplementary Charts for Spatial Error and Dynamic Evolution on BlendedNet}

\begin{figure}[!htbp]
\centering
\includegraphics[width=0.9\textwidth]{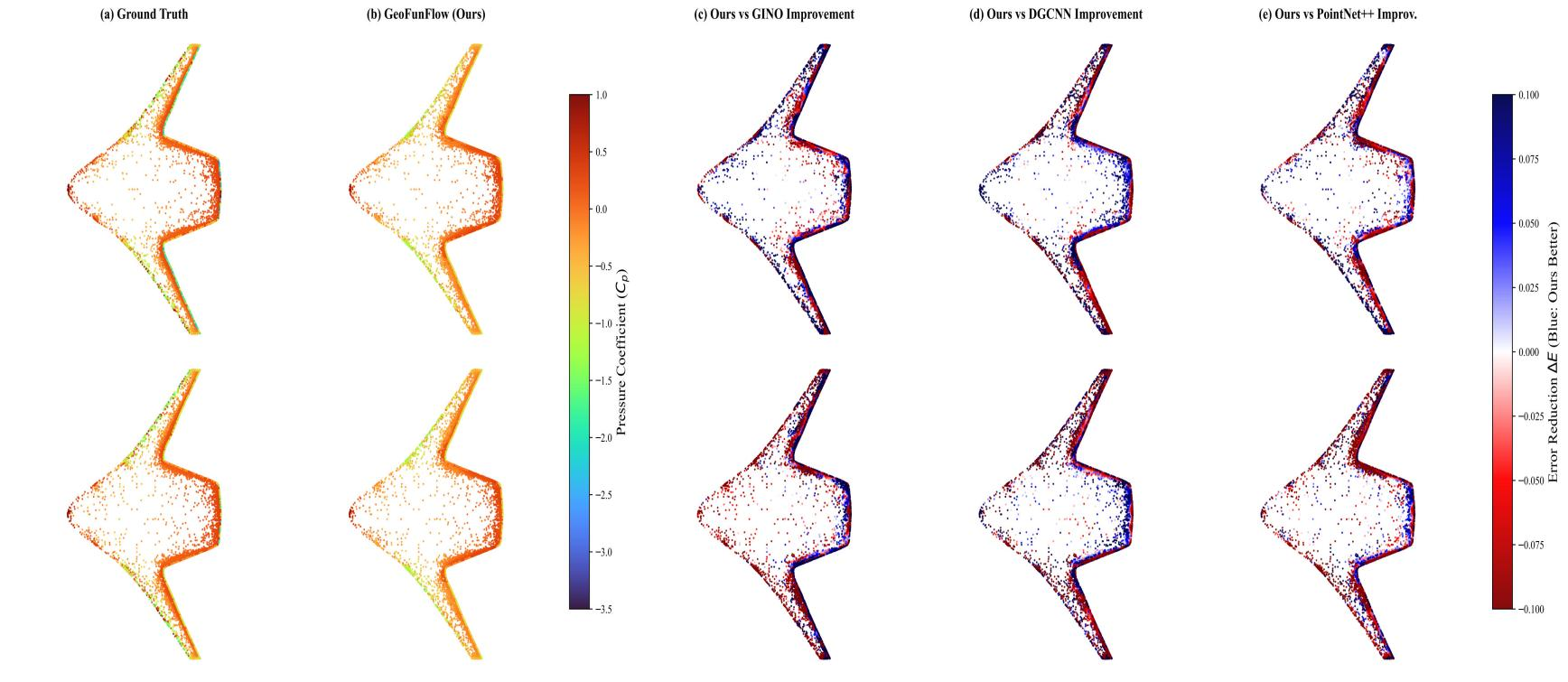} 
\caption{Comparison map of relative error improvement between the GeoFunFlow-3D framework and the baseline models. The blue regions indicate the absolute error reduction of the proposed model in spatial distribution.}
\label{app:error_improvement}
\end{figure}

The spatial accuracy improvement of GeoFunFlow-3D compared to traditional discriminative operators (such as PointNet) needs to be shown intuitively. Therefore, a relative error improvement map is plotted in Figure~\ref{app:error_improvement}. The positive blue regions in the figure indicate that the proposed model has smaller absolute errors at the corresponding nodes. This advantage is clearly observed in regions with geometric curvature mutations.

\begin{figure*}[!htbp]
\centering
\includegraphics[width=0.9\textwidth]{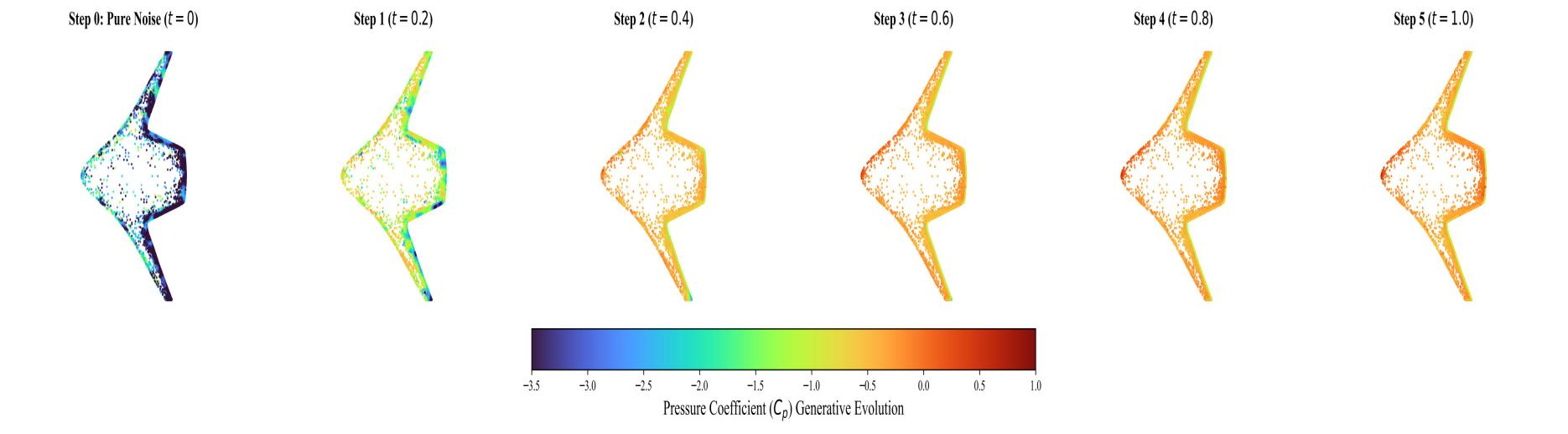} 
\caption{Visualized evolution process of the 3D aerodynamic flow field from pure noise to a physical steady state.}
\label{app:generation}
\end{figure*}

The entire multi-step inverse solving process from Gaussian white noise to a deterministic physical flow field is shown in Figure~\ref{app:generation}. This process is based on the $2$-Wasserstein ($W_2$) geodesic flow matching theory.

\begin{figure*}[!htbp]
\centering
\includegraphics[width=0.9\textwidth]{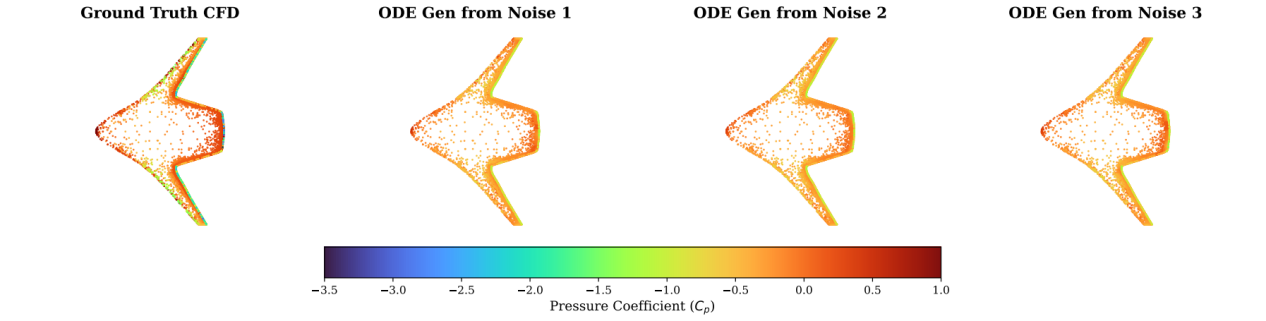} 
\caption{Alignment of the final steady-state generation results using different extreme Gaussian noises as inference starting points.}
\label{app:extreme_noise}
\end{figure*}

The robustness introduced by the dual optimization dynamics was evaluated. To achieve this, completely different extreme Gaussian noises were injected into the latent space. As shown in Figure~\ref{app:extreme_noise}, different initial states can all stably converge to the target physical flow field eventually.

\subsection{Supplementary Charts for Multi-Physics Field Generation on Rotor37}

\begin{figure*}[!htbp]
\centering
\includegraphics[width=0.72\textwidth]{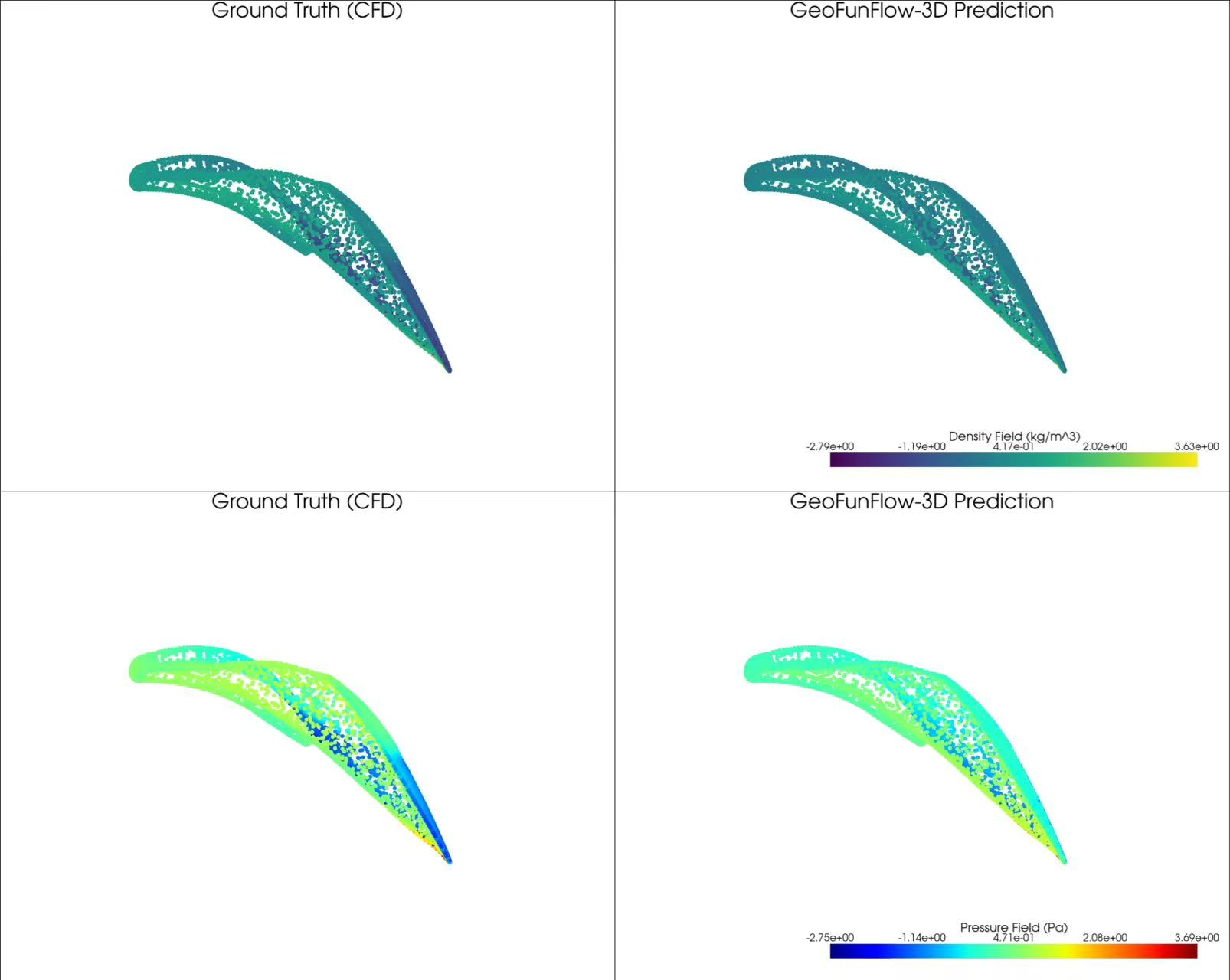} 
\caption{Generation comparison of 3D density and pressure field slices inside the NASA Rotor37 transonic compressor rotor.}
\label{app:rotor_fields}
\end{figure*}

The 3D pressure field was discussed in depth in the main text. Besides this, the density field inside the transonic compressor also exhibits high nonlinearity and transonic shock discontinuities. Joint generation comparison slices of the density and pressure fields are provided in Figure~\ref{app:rotor_fields}. This indicates that the model can effectively suppress local high-frequency numerical oscillations. This is achieved under the protection of total variation (TV) regularization.

\clearpage
\section{Sensitivity Analysis of Tikhonov Regularization Parameter $\epsilon$}
\label{app:tikhonov}

As introduced in Section 3.3.1, the diagonal Tikhonov regularization parameter $\epsilon$ stabilizes the inverse Jacobian multiplier $(R_{i} + \epsilon)^{-1}$. While it ensures numerical stability, we must verify that it does not cause artificial numerical dissipation.

Trade-off and Hyperparameter Selection:
The parameter $\epsilon$ limits the inverse Jacobian $J_{f}^{-1}$ near geometric singularities, such as sharp trailing edges. If $\epsilon$ is too small, gradient calculations become unstable and oscillate. If $\epsilon$ is too large, it acts as numerical viscosity and smooths out physical discontinuities like shocks. In our framework, $\epsilon$ is fixed based on machine precision and the minimum grid spacing, rather than being a trainable parameter. For the BlendedNet and Rotor37 datasets, setting $\epsilon = 10^{-6}$ provides an optimal balance. It prevents $J_{f}^{-1}$ from diverging while remaining negligible compared to the spatial bounding radius $R_{i}$.

Immunity to Numerical Dissipation:
Traditional CFD schemes often add artificial viscosity to the momentum equations. Instead, our regularization only modifies the coordinate mapping operator. Any minor high-frequency smoothing caused by $\epsilon$ is later corrected by the SATO module. As discussed in Section 4.4, SATO generates microscopic basis functions $\mathcal{R}_{SATO}$ to recover sharp gradients. Thus, $\epsilon$ improves stability on complex meshes without reducing the resolution of physical boundaries and shock waves.

\clearpage 
\bibliographystyle{elsarticle-num} 
\bibliography{reference}

\end{document}